\documentclass{article}

\usepackage{amsmath}
\usepackage{amsfonts}
\usepackage{amsthm, amssymb}
\usepackage{eucal}
\usepackage{graphicx}

\newcommand{\colmod}[1]{
\ensuremath{\mathcal{E}_{M}\left(#1\right)}
}

\newcommand{\colneg}[1]{
\ensuremath{\mathcal{N}\left(#1\right)}
}

\newcommand{\col}[1]{
\ensuremath{\mathcal{G}\left(#1\right)}
}

\newcommand{\wf}[1]{
\ensuremath{\rm{WF}(#1)}
}
\newcommand{\ssupp}[1] {
\ensuremath{\rm{sing\ supp}(#1)}
}

\newcommand{\colmap}[2]{
\ensuremath{\mathcal{G}\!\left[#1,#2\right]}
}

\newcommand{\colcomp}[1]{
\ensuremath{\mathcal{G}_c\left(#1\right)}
}

\newcommand{\csub}{
 \subset\!\subset
}

\newcommand{\four}[1]{
\mathcal{F}(#1)
}

\newcommand{\cp}[1]{
\ensuremath{ \rm{CP}\left(#1 \right)}
}

\begin{document}

\title{Microlocal Analysis of generalized pullbacks of Colombeau functions}
\author{Simon Haller, University of Vienna, Faculty of Mathematics}
\date{February 01, 2007}

\maketitle
\begin{abstract}
In distribution theory the pullback of a general distribution by a $C^{\infty}$-function is well-defined
whenever the normal bundle of the $C^{\infty}$-function does not intersect the wavefront set of the distribution.
However, the Colombeau theory of generalized functions allows for a pullback 
by an arbitrary c-bounded generalized function.
It has been shown in previous work that in the case of multiplication of Colombeau functions 
(which is a special case of a $C^{\infty}$ pullback), 
the generalized wave front set of the product satisfies the same inclusion relation as in the distributional case, 
if the factors have their wavefront sets in favorable position. 
We prove a microlocal inclusion relation for the generalized pullback (by a c-bounded generalized map) of Colombeau functions. 
The proof of this result relies on a stationary phase theorem for generalized phase functions, which is given in the Appendix.
Furthermore we study an example (due to Hurd and Sattinger), where the pullback function stems from the generalized characteristic flow 
of a partial differential equation.
\end{abstract}

\newtheorem{definition}{Definition}[section]
\newtheorem{theorem}[definition]{Theorem}
\newtheorem{example}[definition]{Example}
\newtheorem{remark}[definition]{Remark}
\newtheorem{lemma}[definition]{Lemma}
\newtheorem{proposition}[definition]{Proposition}
\newtheorem{diagram}[definition]{Diagram}
\newtheorem{corollary}[definition]{Corollary}

\section{Introduction}
The pullback of a general distribution by a $C^{\infty}$-function in classical distribution theory, 
as defined in \cite[Theorem 8.2.4]{Hoermander:V1}, exists if the normal bundle of the $C^{\infty}$ function
intersected with the wavefront set of the distribution is empty.
These microlocal restrictions reflect also the well-known fact that in general distribution theory one cannot carry out multiplications unrestrictedly,
since the product of two distributions can formally be written as the pullback of a tensor product of the two factors by the diagonal map $\delta:x \mapsto (x,x)$.

Generalized functions in the sense of Colombeau extend distribution theory in a way that it becomes a differential algebra with a
product that preserves the classical product $\cdot: C^{\infty} \times C^{\infty} \rightarrow C^{\infty}$. 
In addition \cite[Proposition 1.2.8]{GKOS:01} states that the Colombeau algebra of generalized functions allows the definition of
a pullback by any c-bounded generalized function.
The classical concept of a wavefront set has been extended to generalized functions of Colombeau type in \cite{DPS2:98,DPS:98,Hoermann:99}.

In \cite{HK:01} the microlocal properties of a product of generalized functions were investigated (which can be interpreted as the pullback of 
a generalized function by a $C^{\infty}$ function). It was shown that the classical microlocal inclusion relations for the wavefront set of a product of distributions (as in \cite[Theorem 8.2.10]{Hoermander:V1}) can be extended to generalized functions, if the wavefront sets of the factors are in favorable position. 

Furthermore the authors provided a counterexample for a product of Colombeau functions with wavefront sets
in unfavorable position that fails the classical microlocal inclusion.

In the present paper we study the wave front set of a generalized pullback of a Colombeau function.

We have divided the paper into four sections:
\begin{itemize}
\item In the first section we introduce the concept of a generalized graph. 
\item In Section 2 we define the transformation of a wave front set by a c-bounded generalized map and the normal bundle of a c-bounded generalized map using the topological concept introduced in the first section.
\item The main theorem is stated in the third section. It is a generalization of \cite[Theorem 8.2.4]{Hoermander:V1}.
\item Section 4 provides two examples: The counterexample mentioned above, from \cite{HK:01}, is investigated in the light of 
our main theorem. We also study the Hurd-Sattinger example (a partial differential equation of first order with non-smooth coefficient) given in \cite{HS:68}, which was further investigated and solved in the Colombeau algebra 
(it is not solvable in $L^1_{\rm{loc}}$, when distributional products are employed) in \cite{HdH:01}. The microlocal properties of the generalized solution can be determined using our main theorem, since the solution is the pullback of the initial data (a generalized function) by the characteristic flow (a c-bounded generalized map).
\item In the Appendix a stationary phase theorem for generalized phase functions is presented, which is crucial for the proof of our main theorem in section 3.
\end{itemize}

\subsection{Notation and basic notions from Colombeau theory}
We use \cite{Colombeau:84,Colombeau:85,GKOS:01,O:92} as the standard references for the foundations of Colombeau theory. In the present paper
we are working with special Colombeau algebras, denoted by $\mathcal{G}^s$ in \cite{GKOS:01}, although here we will drop the superscript 's' to avoid notational overload. 

\subparagraph{The Colombeau Algebra:}
Let us recall the basic construction: A Colomb\-eau 
(generalized) function on some open set $\Omega \subseteq \mathbb{R}^n$ is defined as equivalence class $(u_{\varepsilon})_{\varepsilon}$ of
nets of smooth functions $u_{\varepsilon} \in C^{\infty}(\Omega)$ and $\varepsilon\in]0,1]$ subjected to some asymtotic norm conditions (with respect to $\varepsilon$) for their derivatives on compact sets. We have the following:
\begin{enumerate}
\item \label{labelfirst} Moderate nets $\colmod{\Omega}$: $(u_{\varepsilon})_{\varepsilon}\in C^{\infty}(\Omega)^{]0,1]}$ such that, for all $K\csub\Omega$ and
$\alpha \in \mathbb{N}^n$, there exists $p\in \mathbb{R}$ such that
\begin{equation} \label{modest}
\sup_{x\in K} \|\partial^{\alpha} u_{\varepsilon}(x)\| = O(\varepsilon^{-p})\ \rm{as} \ \varepsilon\rightarrow 0.
\end{equation}
\item Neglible nets $\colneg{\Omega}$: $(u_{\varepsilon})_{\varepsilon}\in C^{\infty}(\Omega)^{]0,1]}$ such that, for all $K\csub\Omega$ and
all $q \in \mathbb{R}$ an estimate
\begin{equation*}
\sup_{x\in K} \|u_{\varepsilon}(x)\| = O(\varepsilon^{q})\ \rm{as} \ \varepsilon\rightarrow 0.
\end{equation*}
holds.
\item \label{labellast} $\colmod{\Omega}$ is a differential algebra with operations defined at fixed $\varepsilon$, $\colneg{\Omega}$ is an ideal
and  $\col{\Omega}:=\colmod{\Omega}/ \colneg{\Omega}$ is the special Colombeau algebra. 
\item If we replace the nets of smooth functions by nets of real numbers in (i)-(iii) we obtain the ring of generalized numbers $\widetilde{\mathbb{R}}$.
\item There are embeddings, $\sigma: C^{\infty}(\Omega) \hookrightarrow \col{\Omega}$ as subalgebra and \\ $\iota: \mathcal{D}'(\Omega) \hookrightarrow \col{\Omega}$ as linear space, commuting with partial derivatives.
\item $\Omega \rightarrow \col{\Omega}$ is a fine sheaf and $\colcomp{\Omega}$ denotes the subalgebra of elements with compact support; using a cut-off in a neighboorhood of the support, one can always obtain representing nets with supports contained in a joint compact set.
\end{enumerate}
\subparagraph{Regular Colombeau functions:} The subalgebra $\mathcal{G}^{\infty}$ of regular Colomb\-eau (generalized) functions 
consists of those elements in $\col{\Omega}$ possessing representatives
such that the estimate (\ref{modest}) holds for a certain $m$ uniformly over all $\alpha \in \mathbb{N}^n$.\\
\subparagraph{Rapdidly decreasing Colombeau functions:}
A Colombeau function in $\col{\mathbb{R}^n}$ is called rapidly decreasing in the directions $\Gamma \subseteq S^{n-1}$, 
if it has a representative with the property that there exists a $N\in\mathbb{N}_0$, such that
for all $p \in \mathbb{N}_0$
\begin{equation*}
\sup_{(\lambda, \xi_1) \in \mathbb{R}^{+} \times \Gamma } (1+\lambda^2)^{p/2}|u_{\varepsilon}(\lambda \xi_1)| = O(\varepsilon^{-N}) \ \rm{as} \ \varepsilon\rightarrow 0
\end{equation*}
holds. 
\subparagraph{Fouriertransform of a compactly supported Colombeau function:}
The Fouriertransform of a $u \in \colcomp{\Omega}$ is a well-defined Colombeau function in $\col{\mathbb{R}^n}$ by 
\begin{equation*}
\mathcal{F}(u) := \left( \int_{\Omega} u_{\varepsilon}(x) e^{-i \langle x, \cdot \rangle} dx\right)_{\varepsilon} + \colneg{\mathbb{R}^n},
\end{equation*}
where $(u_{\varepsilon})_{\varepsilon}$ is a representative of $u$ with joint compact support in $\Omega$.

\subparagraph{$\mathcal{G}^{\infty}$ wavefront set of a Colombeau function:}
If $v \in \colcomp{\Omega}$, we define the set $\Sigma(v) \subset S^{n-1}$ to be the complement of those points having open neighboorhoods $\Gamma\subseteq S^{n-1}$ such that $\four{v}$ is rapidly decreasing in the directions $\Gamma$. $\Sigma(v)$ is a closed subset of $S^{n-1}$.
Now let $u \in \col{\Omega}$. Then we define the cone of irregular directions at $x_0$ by
\begin{eqnarray*}
\Sigma_{x_0}(u) = \bigcap_{\varphi \in C^{\infty}_c(\Omega), \varphi(x_0)\ne 0} \Sigma(\varphi \cdot u ). 
\end{eqnarray*}
Then the (generalized) wave front set of $u$ is the set
\begin{eqnarray*}
\wf{u}:=\left\{ (x,\xi)\in \Omega \times S^{n-1} | \xi \in \Sigma_x(u) \right\}
\end{eqnarray*}
We denote the projection onto the first component by
\begin{eqnarray*}
\ssupp{u}:= \left\{ x\in \Omega | (x,\xi) \in WF(u)\right\}
\end{eqnarray*}
and call this (generalized) singular support of $u$.

\section{The generalized graph}
In the first section we introduce the concept of a generalized graph for a c-bounded generalized map. The generalized graph extends the classical graph of a continous map, in 
the sense that the generalized graph of the embedded map (which is a c-bounded generalized map) coincides with the classical graph. Furthermore the generalized graph is
closed in the product topology of the domain and the image space of the generalized map. This concept will enable us to define the normal bundle of a generalized c-bounded map and the transformation of a wave front set by a generalized c-bounded map (this is carried out in Section \ref{sectionwfset}). Throughout this section we let
 $\Omega_1, \Omega_2$ be open subsets of $\mathbb{R}^n$ resp. $\mathbb{R}^m$.
\begin{definition}
By $\colmap{\Omega_1}{\Omega_2}$ we denote the generalized maps $F\in\mathcal{G}(\Omega_1)^m$ with the property that $F$ is c-bounded on $\Omega_1$ (into $\Omega_2$), i.e. it possesses a representative $(F_{\varepsilon})_{\varepsilon}$ satisfying the condition
\begin{eqnarray*}
\forall K \csub \Omega_1: \exists K' \csub \Omega_2, \exists \varepsilon_0>0 : \textrm{such that } \forall \varepsilon \le \varepsilon_0 : F_{\varepsilon}(K ) \subseteq K'.
\end{eqnarray*}
\end{definition}

\begin{definition}
For $(\lambda_{\varepsilon})_{\varepsilon}$ a net with $\lambda_{\varepsilon} \in K \csub \Omega_2$,  we introduce
\begin{eqnarray*}
\cp{(\lambda_{\varepsilon})_{\varepsilon}}
\end{eqnarray*}  
to denote the set of its clusterpoints.

Let $F\in\colmap{\Omega_1}{\Omega_2}$ be a c-bounded generalized map. Then the set 
\begin{equation*}
\begin{split}
C_{F} := & \left\{(x,y) \in \Omega_1 \times \Omega_2 \mid \exists \text{\ a \ net \ }(x_{\varepsilon})_{\varepsilon} 
\text{\ in\ } \Omega_1: \lim_{\varepsilon\rightarrow 0} x_{\varepsilon}=x \in \Omega_1 \right. \\
&\left. \text{ and } y\in\cp{(F_{\tau(\varepsilon)}(x_{\varepsilon}))_{\varepsilon}}  
 \text{for\ some\  map}\ \right.\\
 &\left. \tau:]0,1] \mapsto ]0,1], \tau(\varepsilon) \le \varepsilon \right\} 
\end{split}
\end{equation*} is called the generalized graph of $F$. Note that we do not consider nets $(x_{\varepsilon})_{\varepsilon}$ in $\Omega_1$ that
converge to the boundary $\partial \Omega_1$.
\end{definition}

\begin{definition}
Let $F\in\colmap{\Omega_1}{\Omega_2}$ be a c-bounded generalized map. 
$F$ is said to be equi-continuous at $x_0$, if
there exists a representative $(F_{\varepsilon})$ which is equi-continuous in $x_0$: For all $\gamma>0$,
there exists some $\delta>0$ such that
\begin{equation*}
|F_{\varepsilon}(x) - F_{\varepsilon}(x_0) | < \gamma
\end{equation*}
holds for all $x \in B_{\delta}(x_0)$ and $\varepsilon \in ]0,1]$.
The generalized map $F$ is called locally equi-continous on some subset $X \subseteq \Omega_1$, if it has a representative which
is equi-continuous in each point $x \in X$.
\end{definition}

\begin{example}
Consider the c-bounded Colombeau function $F$ defined by the representative $(F_{\varepsilon})_{\varepsilon}$ 
\begin{equation*}
F_{\varepsilon} := \sin{(x/\varepsilon)},
\end{equation*}
then it is easy to verify by the definition above, that the generalized graph is
\begin{equation*}
C_F := \mathbb{R} \times [-1,1].
\end{equation*}
\end{example}

\begin{example}
Consider the Colombeau function $F:=\iota(H)$, where $H$ is the Heavyside function. If we set $g:= H \ast \rho$, where $\rho$ is
the mollifier of the embedding $\iota$ with $\int \rho(y) dy=1$, we have that $F_{\varepsilon}(x):=g(x/\varepsilon)$ defines a representative of $F$.
For all nets $(x_{\varepsilon})_{\varepsilon}$ tending to $x_0 \ne 0$, we obtain that $F_{\varepsilon}(x_{\varepsilon}) \rightarrow \rm{sign}(x_0)$ as $\varepsilon \rightarrow 0$.
Let us consider the nets $(x_{\varepsilon})_{\varepsilon}$ tending to $0$: Since we can find zero-nets $(x_{\varepsilon})_{\varepsilon}$ such that $(x_{\varepsilon} / \varepsilon)$ tends to any point in $\mathbb{R} \cup \{\pm \infty\}$ and $g$ is continuous and bounded, we have that
\begin{equation*}
C_F := \{x\in\mathbb{R}\mid x<0\} \times \{0\} \cup \{0\} \times \overline{{\rm Im}(g)} \times \{x\in\mathbb{R}\mid x>0\} \times \{1\}.
\end{equation*}
Note that $\overline{{\rm Im}(g)}$ is a closed interval that contains $[0,1]$.
\end{example}

\begin{proposition}\label{propequi}
Let $F\in\colmap{\Omega_1}{\Omega_2}$ be a c-bounded  generalized map. 
If $F$ is locally equi-continous on ${\Omega_1}$,
then it follows that the generalized graph of $F$ can be determined pointwise by
\begin{equation*}\begin{split}\label{propequistat}
C_F = & \left\{(x,y) \in \Omega_1 \times \Omega_2 \mid y \in \cp{ (F_{\varepsilon}(x))_{\varepsilon}}\ \right\} 
\end{split}
\end{equation*}
\end{proposition}
\begin{proof}
We can choose a locally equi-continous representative $(F_{\varepsilon})_{\varepsilon}$ and proceed along a standard argument:
Let $(x_{\varepsilon})_{\varepsilon},(x'_{\varepsilon})_{\varepsilon}$ be two nets tending to $x_0 \in \Omega_1$. Then we have, using the local equi-continuity of $(F_{\varepsilon})_{\varepsilon}$ on $\Omega_1$, that for any map $\tau:]0,1] \rightarrow ]0,1]$ with $\tau(\varepsilon)<\varepsilon$ and all $\gamma>0$ we can find some $\delta>0$, 
such that the distance is bounded by
\begin{eqnarray*}
|F_{\tau(\varepsilon)}(x_{\varepsilon}) - F_{\tau(\varepsilon)}(x_{\varepsilon}')| < \gamma 
\end{eqnarray*}
for all $x_{\varepsilon}, x'_{\varepsilon} \in B_{\delta/2}(x_0)$. 
Since $\lim_{\varepsilon \rightarrow 0} x_{\varepsilon} =
\lim_{\varepsilon \rightarrow 0} x_{\varepsilon}^{\prime} =x_0$ there exists $\varepsilon' \in ]0,1]$ such that
$|x_{\varepsilon} - x'_{\varepsilon}|<\delta$ holds for all $\varepsilon < \varepsilon'$.
Thus $\lim_{\varepsilon\rightarrow0}F_{\tau(\varepsilon)}(x_{\varepsilon})=\lim_{\varepsilon\rightarrow0}F_{\tau(\varepsilon)}(x_{\varepsilon})$ and it suffices to consider the constant net $x_{\varepsilon}:=x_0$.
For any map $\tau:]0,1] \rightarrow ]0,1]$ with $\tau(\varepsilon)<\varepsilon$ we have that
 $(F_{\tau(\varepsilon)}(x_0))_{\varepsilon}$ is a subnet of $(F_{\varepsilon}(x_0))_{\varepsilon}$, hence (\ref{propequistat}) follows.
\end{proof}

\begin{lemma} \label{nicelemma}
Let $F\in\colmap{\Omega_1}{\Omega_2}$ be a c-bounded generalized map. 
If $(x_0,y_0) \in (\Omega_1 \times \Omega_2) \backslash C_F$, then
there exist neighboorhoods $Y_0 \subseteq \Omega_2$ of $y_0$, $X_0 \subseteq \Omega_1$ of $x_0$ and $\varepsilon'>0$ such that
\begin{equation*}
Y_0 \cap F_{\varepsilon}(X_0) = \emptyset 
\end{equation*}
for all $\varepsilon <\varepsilon'$.
\end{lemma}
\begin{proof}
Let $(x_0,y_0) \in  (\Omega_1 \times \Omega_2) \backslash C_F$. 
The proof proceeds by contradiction: 
Assume that for all neighboords $Y' \subseteq \Omega_2$ of $y_0$ and $X' \subseteq \Omega_1$ of $x_0$ and for
all $\varepsilon'>0$ there exists some $\tau=\tau(X',Y',\varepsilon') < \varepsilon'$ with
\begin{equation*}
Y' \cap F_{\tau(\varepsilon')}(X') \ne \emptyset.
\end{equation*}
Then for all $\varepsilon \in ]0,1]$, by setting $X':=B_{\varepsilon}(x_0) \cap \Omega_1$, $Y':=B_{\varepsilon}(y_0) \cap \Omega_2$ and $\varepsilon':=\varepsilon$,
we can find some $x_{\varepsilon} \in B_{\varepsilon}(x_0) \cap \Omega_1$ and
$\tau(\varepsilon)<\varepsilon$ such that $F_{\tau(\varepsilon)}(x_{\varepsilon}) \in B_{\varepsilon}(y_0) \cap \Omega_2$.
Now $(x_{\varepsilon})_{\varepsilon}$ is a convergent net in $\Omega_1$ with $\lim_{\varepsilon \rightarrow 0} x_{\varepsilon}=x_0 \in \Omega_1$
and the net $(F_{\tau(\varepsilon)}(x_{\varepsilon}))_{\varepsilon}$ converges to $y_0$. Since $\varepsilon \mapsto \tau(\varepsilon)$
defines a map with $\tau(\varepsilon)<\varepsilon$ and $y_0 \in \cp{(F_{\tau(\varepsilon)} (x_{\varepsilon}))_{\varepsilon}}$ we have obtained $(x_0, y_0) \in C_F$, which is
a contradiction.
\end{proof}

\begin{remark} \label{closedproj}
Let $X,Y$ be locally compact topological Hausdorff spaces. 
Assume that $A \subseteq X \times Y$ is a closed set (in the product topology) with the property that
for all compact sets $X_0 \csub X$, it holds that
$A \cap (X_0 \times Y)$ is compact. Then it follows that the projection onto $X$
\begin{equation*}\begin{split}
\pi_X:& X \times Y \rightarrow X \\
& (x,y)\mapsto x
\end{split}
\end{equation*}
has the property that $\pi_X\! \mid_{A}$ is proper and is thus closed. In particular $\pi_X(A)$ is a closed set.
\end{remark}

\begin{proposition}\label{closedproj2}
Let $F\in\colmap{\Omega_1}{\Omega_2}$ be a c-bounded generalized map,
then the generalized graph $C_F$ has the following properties:
\begin{enumerate}
\item It extends the classical notion of a graph (of a continous function) in the following sense: 
If $F$ is the embedding of a continuous function $G := (G_1, ..., G_m) \in C(\Omega_1,\Omega_2)$, i.e. $F=(\iota(G_1),...,\iota(G_m))$,
then the generalized graph of $F$ coincides with the graph of the continous function $G$.
\item
$C_F$ is a closed set in the product topology of $\Omega_1 \times \Omega_2$ and
\item
for any compact set $X_0 \csub \Omega_1$ it holds that  $C_F \cap (X_0 \times \Omega_2)$ is a compact subset of $\Omega_1 \times \Omega_2$. 
\end{enumerate}
\end{proposition}
\begin{proof} 
\begin{trivlist}
\item (i) follows by Proposition \ref{propequi} and the fact that the embedding of a continous function is a locally equi-continous generalized function.
\item (ii) Let $(x_0,y_0) \in (\Omega_1 \times \Omega_2) \backslash C_{F}$, then we know by Lemma \ref{nicelemma} that there exists $X_0 \subseteq \Omega_1, Y_0 \subseteq \Omega_2$
open neighboorhoods of $x_0$ resp. $y_0$ and $\varepsilon'>0$, such that
\begin{equation*}
F_{\varepsilon}(X_0) \cap Y_0 = \emptyset
\end{equation*}
for all $\varepsilon < \varepsilon'$.
Now for any point $x \in X_0$ and any net $(x_{\varepsilon})_{\varepsilon}$ in $\Omega_1$ converging to $x$, there
exists some $\varepsilon'>0$ such that $x_{\varepsilon} \in X_0$ for all $\varepsilon < \varepsilon'$.
It follows that
\begin{equation*}
C_F \cap ( X_0 \times Y_0 ) = \emptyset,
\end{equation*}
so $(\Omega_1 \times \Omega_2) \backslash C_F$ is open in  $\Omega_1 \times \Omega_2$ and hence $C_F$ is closed in the product topology of $\Omega_1 \times \Omega_2$.
\item (iii) Since $F$ is c-bounded we can find some $\varepsilon'>0$ and a compact neighboorhood $K \csub \Omega_1$ of $X_0$ and $K' \csub \Omega_2$,
such that 
\begin{equation*}
F_{\varepsilon}(K) \subseteq K' 
\end{equation*}
for all $\varepsilon < \varepsilon'$. Now if $(x,y) \in C_F \cap (X_0 \times \Omega_2)$, there exists a net $(x_{\varepsilon})_{\varepsilon} $ in $\Omega_1$ converging to $x$, such that $y \in \cp{(F_{\tau(\varepsilon)}(x_{\varepsilon}))_{\varepsilon}}$, where $\tau:]0,1] \rightarrow ]0,1]$ is some
map with $\tau(\varepsilon)<\varepsilon$.
Since $(x_{\varepsilon})_{\varepsilon}$ converges to $x$ and $K$ is a neighboorhood of $x$ there exists $\varepsilon'>0$,
such that $x_{\varepsilon} \in K$ for all $\varepsilon < \varepsilon'$. It follows that $(F_{\tau(\varepsilon)}(x_{\varepsilon}))_{\varepsilon} \in K'$
for all $\varepsilon < \varepsilon'$. So we have that $\cp{(F_{\varepsilon}(x_{\varepsilon}))_{\varepsilon}} \subseteq K'$ and  $(x,y) \in K \times K'$. We conclude that $C_F \cap (X_0 \times \Omega_2) = C_F \cap (X_0 \times \Omega_2) \cap K\times K'$ is a compact subset of $\Omega_1 \times \Omega_2$. 
\end{trivlist}
\end{proof}

\begin{lemma} \label{containlemma}
Let $F\in \colmap{\Omega_1}{\Omega_2}$ be a c-bounded generalized map. Furthermore,
let $X_0 \csub \Omega_1$ and $Y_0 := \pi_2(C_F \cap ( X_0 \times \Omega_2))$.
If $Y$ is a neighboorhood of $Y_0$, then there exists a neighboorhood $X \subseteq \Omega_1$ of $X_0$ and $\varepsilon'>0$, such that
\begin{equation*}
F_{\varepsilon}(X) \subseteq Y
\end{equation*}
for all $\varepsilon < \varepsilon'$. 
\end{lemma}
\begin{proof} 
First we proof the Lemma for $X_0:=\{x_0\}$ containing only a single point $x_0\in\Omega_1$.
We note that $Y_0$ is a compact subset of $\Omega_2$ by Proposition \ref{closedproj2}. 
Assume that the statement of the Lemma does not hold, that is for all neighboorhoods $X'$ of $x_0$ and for
all $\varepsilon>0$ there exists some $\tau(\varepsilon) < \varepsilon$ and $x_{\varepsilon} \in X'$ with
\begin{equation*}
F_{\tau(\varepsilon)}(x_{\varepsilon}) \not \in Y.
\end{equation*}
Then (by setting $X'=B_{\varepsilon}(x_0) \cap \Omega_1$) for all  $\varepsilon \in ]0,1]$, 
we can find a $x_{\varepsilon} \in B_{\varepsilon}(x_0)\cap \Omega_1$ and $\tau(\varepsilon)$ with
$\tau(\varepsilon)<\varepsilon$ such that $ F_{\tau(\varepsilon)}(x_{\varepsilon}) \not \in Y$.

Since $Y$ is a neighboorhood of $Y_0$, we have  that $\overline{Y^c} \cap Y_0 = \emptyset$.
The net $(x_{\varepsilon})_{\varepsilon}$ was choosen to converge to $x_0$ and the net $(F_{\tau(\varepsilon)}(x_{\varepsilon}))_{\varepsilon} \in Y^c$
has the property that the set of clusterpoints $\cp{(F_{\tau(\varepsilon)}(x_{\varepsilon}))_{\varepsilon}}$
is contained in the closure of $Y^c$. 
By the definition of the generalized graph $C_{F}$ we have that $\cp{ (F_{\tau(\varepsilon)}(x_{\varepsilon}))_{\varepsilon} } \subseteq Y_0$,
which contradicts the fact that $\overline{Y^c} \cap Y_0 = \emptyset$. 

Now we consider the general case, when $X_0$ is an arbitrary compact set. Then $Y$ is a neighboorhood of $\pi_2(C_f \cap X_0 \times \Omega_2)$ and
it follows that $Y$ is a neighboorhood of each set $\pi_2(C_f \cap \{ z\} \times \Omega_2)$ for $z\in X_0$. We can apply the first part of the proof for each point $z$ and obtain (open) neighboorhoods $X_z$ of each point $z$ such that $F_{\varepsilon}(X_z) \subseteq Z$ for all $\varepsilon<\varepsilon_z$,  where $\varepsilon_z\in ]0,1]$ depends on $z$.
Since $X_0$ is compact and $(X_z)_{z\in X_0}$ is an open covering we can choose some finite subcovering $(X_{z_k})_{k=1}^l$. Then $X:=\bigcup_{k=1}^l X_{z_k}$
is a neighboorhood of $X_0$ such that $f_{\varepsilon}(X) \subseteq Y$ holds for $\varepsilon < \varepsilon':=\min_{k=1}^l \varepsilon_{z_k}$.
\end{proof}

\begin{corollary} \label{containcorr}
Let $F\in \colmap{\Omega_1}{\Omega_2}$ be a c-bounded generalized map. Furthermore
let $X_0 \csub \Omega_1$ and $Y_0 := \pi_2(C_F \cap X_0\times \Omega_2)$.
If $Y \subseteq \Omega_2$ is some neighboorhood of $Y_0$, then there exists a neighboorhood $X \subseteq \Omega_1$ of $X_0$, such that
\begin{equation}
C_F \cap (X \times \Omega_2) \subseteq X \times Y.
\end{equation} 
\end{corollary}
\begin{proof}
If $X'$ is the neighboorhood constructed in Lemma \ref{containlemma}, then we can find some $\delta>0$, such 
that $B_{\delta}(X_0) \subseteq X'$. Set $X:= B_{\delta/2}(X_0)$, then 
the statement is true, since $X'$ is a neighboorhood for each point of $X$.
\end{proof}

\section{Transformation of wave-front sets}\label{sectionwfset}

In this section we consider a c-bounded generalized map $f\in \colmap{\Omega_1}{\Omega_2}$, where $\Omega_1 \subseteq \mathbb{R}^n,\ \Omega_2 \subseteq \mathbb{R}^m$ are open sets and $\Gamma \subseteq \Omega_2 \times S^{m-1}$ is a closed set. We define the normal bundle of a c-bounded generalized map and consider the transformation of wavefront sets by c-bounded generalized maps.

\begin{definition}
Let $y\in \Omega_2$ and $Y \csub \Omega_2$ be some compact set. Then we use the notation
\begin{equation*}
\Gamma_y := \{\eta \in S^{m-1} \mid (y,\eta) \in \Gamma \}
\end{equation*}
and
\begin{equation*}
\Gamma_Y := \{\eta \in S^{m-1} \mid \exists y\in Y: (y,\eta) \in \Gamma \}=\bigcup_{y\in Y} \Gamma_y.
\end{equation*}
Note that both sets are closed subsets of $S^{m-1}$, since
\begin{equation*}
\Gamma_Y= \pi_{2}(\Gamma \cap (Y \times S^{m-1}))
\end{equation*}
is the image of a compact set under the continous projection \\$\pi_2: (y,\eta) \mapsto \eta$.
 \end{definition}

\begin{definition}
According to \cite[Definition 2.2]{HO:04} a net $(\lambda_{\varepsilon})_{\varepsilon}$ in $\mathbb{R}$  is said to be of slow-scale, if
\begin{eqnarray*}
\exists \varepsilon'\in]0,1]: \forall t \ge 0, \exists C_t>0 \rm{\ such\ that\ } |\lambda_{\varepsilon}|^t \le C_t \varepsilon^{-1}, \rm{\ for\ all\ } \varepsilon <  \varepsilon'
\end{eqnarray*}
holds.
\end{definition}

\begin{lemma}\label{lemmacone3}Let $\Omega \subseteq \mathbb{R}^m$ be an open set and let $\Gamma$ be some closed subset of $\Omega \times S^{m-1}$. 
Furthermore let $Y_0 \csub \Omega$ be a compact set. If $V \subseteq S^{m-1}$ is some closed neighborhood of $\Gamma_{Y_0}$, then 
there exists a neighboorhood $Y$ of $Y_0$ such that 
\begin{equation*}
\Gamma_{Y} \subseteq V
\end{equation*}
holds.
\end{lemma}

\begin{proof}
We note that $\Gamma_{Y_0}$ is a compact subset of $S^{m-1}$.
For any compact set $L\csub \mathbb{R}^m$  and $\lambda>0$ we define the set $K_{\lambda}(L):=
\{z\in \mathbb{R}^m \mid \exists z' \in L: |z-z'|\le \lambda \}$.
Since $V$ is a compact neighboorhood of the compact set $\Gamma_{Y_0}$ there exists some $\delta >0$
such that $K_{\delta}(\Gamma_{Y_0}) \cap S^{m-1} \cap (S^{m-1} \backslash V) =  \emptyset$.  We proove the statement by contradiction:
Assume the statement is false, then we have that
for all $\varepsilon\in ]0,1]$ there exists a point $y_{\varepsilon} \in K_{\delta \cdot \varepsilon}(Y_0)$ with the property that $\Gamma_{y_{\varepsilon}} \not \subseteq V$. Therefore we can find some $\xi_{\varepsilon} \in \Gamma_{y_{\varepsilon}} \backslash V \ne \emptyset$
for each $\varepsilon\in]0,1]$. This leads to a net 
$(y_{\varepsilon}, \xi_{\varepsilon})_{\varepsilon}$ contained in the compact set $\Gamma\cap (S_{\delta}(Y_0) \times S^{m-1})$.
Thus the net $(y_\varepsilon, \xi_{\varepsilon})_{\varepsilon}$ has a finer net  $(y_{\tau(\varepsilon)},\xi_{\tau(\varepsilon)})_{\varepsilon}$ that converges to some $(y_0,\xi_0) \in \Gamma \cap (K_{\delta}(Y_0) \times S^{m-1})$ and $\xi_0 \in \Gamma_{y_0}$. Since $y_{\varepsilon} \in K_{\varepsilon \delta}(Y_0)$, it follows that the clusterpoint $y_0$ is contained in $Y_0$. We conclude that $\xi_0 \in \Gamma_{y_0} \subseteq \Gamma_{Y_0}$. We have that there exists $\varepsilon_1\in]0,1]$ such that
\begin{eqnarray*}
\max{(|\xi_{\tau(\varepsilon)} -\xi_0|,|y_{\tau(\varepsilon)}-y_0|)} <\delta
\end{eqnarray*}
holds for all $\varepsilon < \varepsilon_1$. We have choosen $\delta$ such that $\{\xi \in S^{m-1} \mid \exists \xi_0\in \Gamma_{Y_0}: |\xi-\xi_0| < \delta \} \subseteq V$. It follows that $\xi_{\tau(\varepsilon)} \in V$ for all $\varepsilon<\varepsilon_1$, which contradicts the choice
$\xi_{\varepsilon} \in \Gamma_{y_{\varepsilon}} \backslash V$.
\end{proof}

Although the pullback of a Colombeau function $u$ by any c-bounded generalized map $f$ is well-defined, we cannot derive a general microlocal inlusion for the pullback without requiring further properties for the generalized map $f$.
We define an open subdomain $D_f$ of $\Omega_1 \times S^{m-1}$, where the generalized map $(x,\eta) \mapsto {}^T\! df_{\varepsilon}(x) \eta$ has certain properties which are needed to obtain a microlocal inclusion relation. This leads to the notion of a generalized normal bundle.

\begin{definition} \label{unfavdef}
Let $f\in\colmap{\Omega_1}{\Omega_2}$ be a c-bounded generalized map, then 
we define the open set $D_f$ by 
\begin{equation*} \begin{split}
D_f := &\{ (x,\eta) \in \Omega_1 \times S^{m-1} \mid \exists\  \rm{\ neighboorhood\ }X \times V \subseteq \Omega_1 \times S^{m-1} \rm{\ of\ } (x,\eta) \\
&\rm{\ and\ a\ positive\ net\ of\ slow\ growth\ } 
(\sigma_{\varepsilon})_{\varepsilon}, \exists \varepsilon'\in ]0,1]:\\
&\inf_{(x,\eta)\in X\times V} 
   |\sigma_{\varepsilon}  {}^T\! df_{\varepsilon}(x) \eta | \ge 1  \ \rm{and}\\
&\sup_{(x,\eta) \in X\times V^{\perp}} |\sigma_{\varepsilon}  {}^T\! df_{\varepsilon}(x) \eta | \le 1 \rm{\ for\ all\ } \varepsilon < \varepsilon'  \},
\end{split}
\end{equation*}
where $V^{\perp}:= \{\eta \in S^{m-1} \mid \exists \eta_0 \in V: \langle \eta, \eta_0 \rangle = 0 \}$. 
Then the generalized normal bundle of $f$ is defined by
\begin{equation*}
N_f := \{(y,\eta) \in \Omega_2 \times S^{m-1} \mid (x,y) \in C_{f}, (x,\eta) \not \in D_f \}.
\end{equation*}
This is the analog to the classical normal bundle.
Furthermore, we define the wave-front unfavorable support of $f$ with respect to a closed set $\Gamma\subseteq \Omega_2 \times S^{m-1}$ by
\begin{eqnarray*} 
U_{f}(\Gamma):= \{ x\in \mathbb{R}^n \mid  (x,y) \in C_{f}, (x,\eta) \not \in D_f, (y,\eta) \in \Gamma \}. 
\end{eqnarray*}
\end{definition}

\begin{example}\label{mgexample}
We consider the c-bounded generalized map defined by 
\begin{eqnarray*}
f_{\varepsilon}(x,y)=(x+ \gamma_{\varepsilon} y, x-\gamma_{\varepsilon} y),
\end{eqnarray*}
where $\lim_{\varepsilon \rightarrow 0} \gamma_{\varepsilon} = 0$.
\end{example}
Then the transposed Jacobian is
\begin{eqnarray*}
{}^T\!df_{\varepsilon}(x,y) := \left( \begin{array}{cc} 
1 &  1 \\
\gamma_{\varepsilon} & -\gamma_{\varepsilon}
\end{array}\right)
\end{eqnarray*}
and ${}^T\!df_{\varepsilon}(x,y) \eta= (\eta_1 +\eta_2, \gamma_{\varepsilon}(\eta_1-\eta_2))$.
Let $\tilde{\eta}:=(\pm 1, \mp 1)$ and $\tilde{\eta}^{\perp}:=(\pm 1,\pm 1)$ then
\begin{eqnarray*}
a_{\varepsilon}(x,\tilde{\eta})&:=&| \sigma_{\varepsilon} {}^T\!df_{\varepsilon}(x) \tilde{\eta}| = \gamma_{\varepsilon} \sigma_{\varepsilon} 2\\
b_{\varepsilon}(x,\tilde{\eta})&:=&| \sigma_{\varepsilon} {}^T\!df_{\varepsilon}(x) \tilde{\eta}^{\perp}| = \sigma_{\varepsilon} 2
\end{eqnarray*}
shows that for any choice of a slow-scaling net $\sigma_{\varepsilon}$  either $a_{\varepsilon} \rightarrow 0, b_{\varepsilon} \rightarrow 1 \in \mathbb{R}^+$ or
$a_{\varepsilon} \rightarrow 1 \in\mathbb{R}^+, b_{\varepsilon} \rightarrow \infty$ (the second case is ruled out if $\gamma_{\varepsilon}$ is not a slow-scaled net).
For all other $\eta \not \in S^{n-1} \backslash (\pm 1 ,\mp 1)$ it is easy to find a neighboorhood $V$ of $\eta$ and
some slow-scaling net $\sigma_{\varepsilon}$,
such that $\inf_{\mathbb{R}^2\times V} a_{\varepsilon}(x,\eta) \ge 1$ and 
$\sup_{\mathbb{R}^2\times V^{\perp}} a_{\varepsilon}(x,\eta) \le 1$ for $\varepsilon$ small enough. 
It follows that $D_f:= \mathbb{R}^2 \times S^1 \backslash (\pm 1, \mp1)$.

\begin{lemma}
The generalized normal bundle $N_f$ and the wave-front unfavorable support $U_f(\Gamma)$ of $f$ (with respect to $\Gamma$) are closed sets.
If $N_f \cap \Gamma =\emptyset$, then $U_f(\Gamma) = \emptyset$.
\end{lemma}
\begin{proof}
Let $x_0 \not \in U_{f}(\Gamma)$, then since $D_f$ is open there exist neighboorhoods 
$X_0$ resp. $V_0$ of $x_0$ resp. $\Gamma_{Y_0}$ such that $X_0 \times V_0 \subseteq D_f$.
By Lemma \ref{lemmacone3} there exists a neighboorhood $Y$ of $Y_0$, such that
$V_0$ is still a neighboorhood of $\Gamma_Y$. Corollary \ref{containcorr} provides a neighboorhood $X'$ of
$x_0$ such that $C_f \cap (X' \times \Omega_2) \subseteq X' \times Y$.
Now we can choose a smaller neighboorhood $X_1 \subset X' \cap X_0$ (such that $X' \cap X_0$ is a neighboorhood of $X_1$) of $x_0$, 
such that $X_1 \times V_0 \subseteq D_f$ and $C_F \cap (X_1 \times \Omega_2) \subseteq X_1 \times Y$.  
For all $x_1\in X_1$ and $y\in \Omega_2$ with $(x_1,y)\in C_f$, it holds that $(x_1,y) \in X_1 \times Y$ and if $(y,\eta) \in \Gamma$ it follows
that $\eta \in \Gamma_Y \subseteq V$ and thus $(x_1,\eta) \in D_f$. It follows that $x_1\not\in U_f$ for all $x_1\in X_1$, $X_1 \cap U_f =\emptyset$ and thus $U_f(\Gamma)^c$ is an open set.
\end{proof}

A substantial step in the proof of the main theorem applies a generalized stationary phase theorem (in the Appendix). In order to obtain a lower bound for the gradient of the occuring phase function we have to consider the map
\begin{eqnarray*}
\Omega_1 \times S^{m-1} \times S^{n-1} & \rightarrow& S^{n-1} \\
(x,\eta,\xi)	&	\mapsto & \left|\frac{{}^T df_{\varepsilon}(x) \eta}{|{}^T df_{\varepsilon}(x) \eta|} -\xi \right|.
\end{eqnarray*}
So we introduce the following notation:

\begin{lemma}\label{notionlemma}
Consider
\begin{eqnarray*}
M_{\varepsilon}(x,\eta):=\frac{{}^T\!df_{\varepsilon}(x) \eta }{|{}^T\!df_{\varepsilon}(x) \eta|}
\end{eqnarray*}
on the domain $D_f$, then $\eta \mapsto M_{\varepsilon}(x,\eta)$ defines an equi-continous Colombeau function in the $\eta$ variable at fixed $x$. 
\end{lemma}
\begin{proof}
The map $g:\xi \rightarrow \frac{\xi}{|\xi|}$
is equi-continous on $\mathbb{R}^n \backslash B_{\delta}(0)$ for any fixed $\delta>0$ since
\begin{eqnarray*}
\left| \frac{\xi}{|\xi|} - \frac{\eta}{|\eta|} \right| = \frac{1}{|\xi|} \left| \xi - \frac{\eta|\xi|}{|\eta|} \right| =\frac{1}{|\xi|}
\left| \xi -\eta + \frac{\eta}{|\eta|} \left( |\eta|-|\xi| \right) \right| \le \frac{2}{|\xi|} |\xi-\eta|.
\end{eqnarray*}
Let $(x_0,\eta_0) \in D_{f}$, then we can find some neighboorhood $X\times V \subseteq \Omega_1 \times S^{m-1}$ of $(x_0,\eta_0)$ such that
\begin{eqnarray*}
\inf_{(x,\eta)\in X\times V}  | \sigma_{\varepsilon}  {}^T\! df_{\varepsilon}(x) \eta|  &\ge& 1 \rm{\ and\ }\sup_{(x,\eta) \in X\times V^{\perp}} |  \sigma_{\varepsilon}  {}^T\! df_{\varepsilon}(x) \eta | \le 1
\end{eqnarray*}
for some slow-scaling net $(\sigma_{\varepsilon})_{\varepsilon}$.
Then we conlude that
\begin{equation*}\begin{split}
&|M_{\varepsilon}(x,\eta)- M_{\varepsilon}(x,\xi)| \le \frac{2}{  | {}^T df_{\varepsilon}(x) \eta |} \left({}^T\! df_{\varepsilon}(x) \eta-{}^T\! df_{\varepsilon}(x) \xi  \right) \\ &=  \frac{2}{ |  \sigma_{\varepsilon} {}^T\! df_{\varepsilon}(x) \eta| } 
\left| \sigma_{\varepsilon} {}^T\! df_{\varepsilon}(x) (\eta-\xi)  \right| \le
\frac{ \sup_{(x,\zeta) \in X\times V^{\perp}} |\sigma_{\varepsilon} {}^T\! df_{\varepsilon}(x) \zeta |}{ 
 \inf_{(x,\eta)\in X\times V}   |\sigma_{\varepsilon} {}^T\! df_{\varepsilon}(x) \eta| } \\
&\cdot |\xi -\eta| \le  2 |\xi- \eta|
\end{split}
\end{equation*}
uniformly for all $x\in X_0$ and $\xi,\eta \in V$.
\end{proof}

\begin{remark}
From the preceding Lemma we obtain that if $\eta_{\varepsilon} \rightarrow \eta_0$ and $x_{\varepsilon} \rightarrow x_0$, then
\begin{eqnarray*}
\cp{ (M_{\varepsilon}(x_{\varepsilon},\eta_{\varepsilon} ))_{\varepsilon}} = \cp{ (M_{\varepsilon}(x_{\varepsilon}, \eta_0 ))_{\varepsilon}},
\end{eqnarray*}which simplifies the determination of the generalized graph $C_{M}$.
\end{remark}

\begin{definition}
Let $f\in\colmap{\Omega_1}{\Omega_2}$ be a c-bounded generalized map and assume that $D_f \ne \emptyset$.
We define the transformation map $\Phi_f$ by 
\begin{equation*} \begin{split}
\Phi_f:& D_f \rightarrow \Omega_2 \times S^{n-1} \\
& (x,\eta) \mapsto \left( f_{\varepsilon}(x), \frac{{}^T\!df_{\varepsilon}(x) \eta }{|{}^T\!df_{\varepsilon}(x) \eta|} \right).
\end{split}
\end{equation*}
Consider the projections
\begin{equation*}\begin{split}
\Pi_{\Phi_f,1}: & C_{\phi_f} \rightarrow \Omega_1 \times S^{n-1},\ (x, \eta ,y, \xi) \mapsto (x,\xi)\\
\Pi_{\Phi_f,2} : & C_{\phi_f} \rightarrow \Omega_2 \times S^{m-1},\ (x, \eta ,y, \xi) \mapsto (y,\eta)\\
\end{split}
\end{equation*}
and define the set $f^{\ast} \Gamma$ by
\begin{equation*}\begin{split}
f^{\ast} \Gamma &:= \Pi_{\Phi_f,1} \circ \Pi_{\Phi_f,2}^{-1}(\Gamma) \\
&= \{ (x,\xi) \in \Omega_1 \times S^{n-1} | (x,\eta,y,\xi ) \in C_{\Phi_f}, (y,\eta) \in \Gamma \}.
\end{split}
\end{equation*}
In the case of $D_f = \emptyset$, we set $f^{\ast} \Gamma := \Omega_1 \times S^{n-1}$.

\begin{lemma}
Let $f\in\colmap{\Omega_1}{\Omega_2}$ be a c-bounded generalized map and let $\Gamma \subset \Omega_2 \times S^{m-1}$ be a closed set in the product topology of $\Omega_2 \times S^{m-1}$. Then the set $f^{\ast}\Gamma$ is closed in the product topology of $\Omega_1 \times S^{n-1}$. 
\end{lemma}
\end{definition}
\begin{proof}
The maps
\begin{equation*}\begin{split}
\Pi_{1}: \Omega_1 \times S^{m-1} \times \Omega_2 \times S^{n-1} &\rightarrow \Omega_1 \times S^{n-1},\\ (x,\eta,y,\xi)&\mapsto (x,\xi) \\
\Pi_{2}: \Omega_1 \times S^{m-1} \times \Omega_2 \times S^{n-1} &\rightarrow \Omega_2 \times S^{m-1},\\ (x,\eta,y,\xi)&\mapsto (y,\eta)
\end{split}
\end{equation*}
are continous projections and we have that 
\begin{equation*}
\Pi_{\Phi_f,1} = \Pi_1\mid_{C_{\Phi_f}} 
\end{equation*} and
\begin{equation*}
\Pi_{\Phi_f,2} = \Pi_2\mid_{C_{\Phi_f}}.
\end{equation*}
$C_{ \Phi_f }$ has the following property:
For all $K \csub \Omega_1 \times S^{n-1}$ there exists some $K' \csub \Omega_1 \times S^{m-1}
\times \Omega_2 \times S^{n-1}$ such that $C_{\Phi_f} \cap \Pi_1^{-1}(K) \subseteq K'$ (cf. Proposition \ref{closedproj2}).
By Remark \ref{closedproj} the map $\Pi_{\Phi_f,1}$ is thus closed. It follows that $\Pi_{\Phi_f,1} \circ \Pi_{\Phi_f,2}^{-1}(\Gamma)$ is a closed set in the subset topology of $\Omega_1 \times S^{n-1}$.
\end{proof}

\section{Generalized pullbacks of Colombeau functions}
In this section we prove the main result which gives a microlocal inclusion relation for the generalized pullback of a Colombeau function.
The proof of the theorem relies on a generalized stationary phase theorem, the details of which are discussed in the Appendix.

\begin{definition}
Let $f \in \colmap{\Omega_1}{\Omega_2}$ be a c-bounded generalized map. Then we call $f$ slow-scaled in all
derivatives on the open set $X_0\subseteq \Omega_1$, if for all $\alpha\in\mathbb{N}^n$ there exists slow-scaled nets $(r_{\alpha,\varepsilon})_{\varepsilon}$
such that
\begin{eqnarray}\label{slowscale}
\sup_{X_0} |\partial^{\alpha} f_{\varepsilon}(x)| \le C_{\alpha} r_{\alpha,\varepsilon} \ \rm{as} \ \varepsilon \rightarrow 0
\end{eqnarray}
holds, where $C_{\alpha}$ are constants. Furthermore we call
$f$ slow-scaled in all derivatives at $x_0\in\Omega_1$, if there exists a neighboorhood $X_0$ of $x_0$ such that (\ref{slowscale}) holds for all $\alpha \in \mathbb{N}^n$.

Define the sets
\begin{eqnarray*}
S_f := \{ x \in \Omega_1 \mid f \rm{\ is\ slow-scaled\ in\ all\ derivatives\ at\ } x\}
\end{eqnarray*}
and 
\begin{eqnarray*}
K_f(u) := \bigcap_{k \in {{\widetilde{\mathbb{R}}}^m } } \pi_1(C_f \cap \rm{supp}(u-k) \times \Omega_2).
\end{eqnarray*}
\end{definition}

\begin{remark}
The sets $\pi_1(C_f \cap \Omega_1 \times \rm{supp}(u-k))$ are closed because $\pi_1\! \mid_{C_f}$ is a proper map and the sets 
$\rm{supp}(u-k)\subset \Omega_1$ and $C_f \subseteq \Omega_1 \times \Omega_2$ are closed. Thus the 
set $K_f(u)$ is closed in the relative topology of $\Omega_1 \subseteq \mathbb{R}^n$.
\end{remark}

\begin{lemma}\label{lemmacool}
Let $f \in \colmap{\Omega_1}{\Omega_2}$ be a c-bounded generalized map. Furthermore let $x_0 \not \in U_f(\Gamma) \cup S_f^c$ and $Y_0 := \pi_2(C_f \cap \{x_0\}\times \Omega_2)$. If $W$ is some open neighboorhood of $(f^{\ast} \Gamma)_{x_0}$ , 
then there exist neighboorhoods $X$ of $x_0$, $V$ of $\Gamma_{Y_0}$ 
and $Y$ of $Y_0$ 
with the following properties:
\begin{equation*}\begin{split}
&X \times V \subseteq D_f(\Gamma),\ (f^{\ast} \Gamma)_X \subseteq W,\
f_{\varepsilon}(X) \subseteq Y \ \rm{for\ all\ }\ \varepsilon < \varepsilon', \\ 
&(f_{\varepsilon}) \rm{\ is\ slow-scaling\ in\ all\ derivatives\ on\ } X \rm{\ and\ } 
\Gamma_{Y} \subseteq V \end{split}
\end{equation*}
Furthermore there exists a positive constant $c>0$, such that
\begin{equation*}
\inf_{(x,\eta,\xi)\in X \times V \times W^c} \left|M_{\varepsilon}(x,\eta) - \xi\right|>c \ \rm{for\ all\ }\varepsilon< \varepsilon'
\end{equation*}
holds (we are using the notation from Lemma \ref{notionlemma}).
\end{lemma}
\begin{proof}
From $x_0 \not \in U_f(\Gamma)$ it follows that $\{x_0 \}\times \Gamma_{Y_0}$
is a compact subset of the open set $D_f$. Thus we can find some
neighboorhood $X' \times V' \subseteq D_f$ of $\{x_0\} \times \Gamma_{Y_0}$.
Since $x_0 \in S_f$ we can assume without loss of generality that $f$ is slow-scale in all derivatives on
the compact set $X'$.

Since $W$ is a neighboorhood of $(f^{\ast}\Gamma)_{x_0}$ we have that
\begin{eqnarray*}
(f^{\ast}\Gamma)_{x_0}  = \pi_2( C_{M} \cap \{x_0\} \times \Gamma_{Y_0} \times S^{n-1}) \csub W
\end{eqnarray*}
where $C_M$ denotes the generalized graph of the generalized map defined by $(x,\eta) \mapsto M_{\varepsilon}(x,\eta)$ on the open domain $D_f$.
By Lemma \ref{containlemma}
there exist neighboorhoods $X'', V''$ of $x_0$ resp. $\Gamma_{Y_0}$ such that
\begin{eqnarray} \label{propcool}
M_{\varepsilon}(X'' \times V'') \subseteq W 
\end{eqnarray}
for all $\varepsilon < \varepsilon'$. Let $V:= V' \cap V''$, then by Lemma \ref{lemmacone3} we can find some neighboorhood $Y$ of $Y_0$ such that
$\Gamma_Y \subseteq V$. By Lemma \ref{containlemma} there exists some neighboorhood $X'''$ of $x_0$ such 
that $f_{\varepsilon}(X''') \subseteq Y$ for small $\varepsilon$. Let $X:=X' \cap X'' \cap X'''$ and
$Z:= X \times V \times W^c$. Then there exists $(x_{\varepsilon},y_{\varepsilon},\xi_{\varepsilon})\in Z$ (note that $(x,\eta,\xi) \mapsto |M_{\varepsilon}(x,\eta) - \xi |$ is a continuous function for each $\varepsilon \in ]0,1]$ and $Z$ is a compact set) such that
\begin{equation*}
c_{\varepsilon} := \inf_{(x,\eta,\xi) \in Z} |M_{\varepsilon}(x,\eta) - \xi | = |M_{\varepsilon}(x_{\varepsilon},\eta_{\varepsilon}) - \xi_{\varepsilon} |
\end{equation*}
holds for some net $(x_{\varepsilon}, \eta_{\varepsilon}, \xi_{\varepsilon})_{\varepsilon}$.
By (\ref{propcool}) we have that $c_{\varepsilon} > c >0$ holds for all $\varepsilon < \varepsilon'$, where $c$ is a constant.
\end{proof}

\begin{theorem} \label{maintheorem}
Let $f \in \colmap{\Omega_1}{\Omega_2}$ be a c-bounded generalized map with representative $(f_{\varepsilon})_{\varepsilon}$.
For $u \in \col{ \Omega_2}$ with representative $(u_{\varepsilon})_{\varepsilon}$, we define the pullback  $f^{\ast} u$  by
\begin{eqnarray*}
f^{\ast} u := (u_{\varepsilon}(f_{\varepsilon}(x)))_{\varepsilon} +\mathcal{N}(\Omega_1).
\end{eqnarray*}
Then $f^{\ast}u$ is well-defined and the microlocal inclusion relation
\begin{eqnarray*}
\wf{f^{\ast} u} \subseteq f^{\ast} \wf{u}\ \bigcup\ \left(U_f(\wf{u}) \times S^{n-1}\right)\ \bigcup\ \left(\left(K_f(u) \cap S_f^c \right) \times S^{n-1}\right)
\end{eqnarray*}
holds.
\end{theorem}
\begin{proof}
According to \cite[Proposition 1.2.8]{GKOS:01} the pullback $f^{\ast} u$ is a well-defined Colombeau function in $\mathcal{G}(\Omega_1)$.

Let $x_0 \in \Omega_1$ and $Y_0 =\{ y \in \Omega_2 \mid (x_0, y) \in C_f(\Omega_1) \}$ and
$\xi_0$.
We set $\Gamma:=\wf{ u }$ and $\Gamma_y:=\{\eta \mid (y,\eta) \in \Gamma\}$.\\
Suppose that $(x_0,\xi_0) \not \in f^{\ast} \Gamma  \bigcup U_f(\Gamma) \times S^{n-1} \bigcup \left(K_f(u) \cap S_f^c \right) \times S^{n-1}$, so 
$\xi_0 \not \in (f^{\ast} \Gamma)_{x_0}$ and $x_0 \not \in U_f(\Gamma) \bigcup K_f(u) \cap S_f^c$.
We distinguish the cases $x_0 \not\in K_f(u) \cup U_f(\Gamma)$ and $x_0 \not\in S_f^c \cup U_f(\Gamma)$:

In the first case we have that there exists a generalized constant $k\in \widetilde{\mathbb{R}}^n$ with representative $(k_{\varepsilon})_{\varepsilon}$, 
such that $ C_f \cap \{x_0\}\times \rm{supp}(u-k) = \emptyset$. Hence
for all $y \in Y_0$ it holds that $y \not \in \rm{supp}(u-k)$.
Since $Y_0$ and $\rm{supp}(u-k)$ are two disjoint closed sets,
we can find a closed neighboorhood $Y$ of $Y_0$ such that $Y \cap \rm{supp}(u)=\emptyset$.
Let $\chi \in C_{c}^{\infty}(\Omega_2)$ with the property that $\chi \equiv 1$ on some compact neighboorhood $Y' \subseteq Y^{\circ}$ of $Y_0$
and $\rm{supp}(\chi) \subseteq Y$. By Lemma \ref{containlemma}
there exists a neighboorhood $X'$ of $x_0$ and $\varepsilon'\in]0,1]$ such that $f_{\varepsilon}(X') \subseteq Y'$ for $\varepsilon<\varepsilon'$.

In order to show that $\xi_0 \not \in (\wf{f^{\ast}u})_{x_0}$ for $x_0 \not\in K_f(u)$ we have to find a smooth function $\varphi$ with support on a neighboorhood of $x_0$, such that $\mathcal{F}(f^{\ast}u \varphi)$ is rapidly decreasing on some neighboorhood of $\xi_0$.

We choose $\varphi$ to be a smooth function with $\rm{supp}(\varphi) \subseteq X'$. For all $\varepsilon < \varepsilon'$ the
identity $f_{\varepsilon}^{\ast} u \cdot \varphi = f_{\varepsilon}^{\ast}(\chi u) \cdot \varphi$ holds, since the functions $\chi$ and $\varphi$ where choosen
such that $ (\chi \circ f_{\varepsilon})\cdot \varphi \equiv 1$ for $\varepsilon < \varepsilon'$.
We have that $\chi \cdot (u_{\varepsilon}-k_{\varepsilon}) \in \colneg{\Omega_2}$ and it follows that 
$f^{\ast}u\cdot \varphi$ is $\mathcal{G}^{\infty}$, so $\xi_0 \not \in (\wf{f^{\ast}u})_{x_0}$.\\

Now we consider the case $x_0 \not \in U_f(\Gamma) \cup S_f^c$.
Let $W$ be an open set in $S^{n-1}$ such that $W^c$ is a neighboorhood of $\xi_0$.
Lemma \ref{lemmacool} implies that there exist neighboorhoods $X$ of $x_0$, $Y$ of $Y_0$
and $V$ of $\Gamma_{Y_0}$ such that
$(f^{\ast} \Gamma)_X \subseteq W$, $f_{\varepsilon}(X) \subseteq Y \ \rm{for\ all\ }\ \varepsilon < \varepsilon'$, 
$\Gamma_{Y} \subseteq V$,
$(f_{\varepsilon})$ is slow-scaling in all derivatives on $X$ and 
$\inf_{(x,\eta_1,\xi_1) \in X\times V \times W^c}  \left|M_{\varepsilon}(x,\eta_1) -  \xi_1 \right| > d >0$
holds for all $\varepsilon < \varepsilon'$, where $d$ is some positive constant.\\
Let $\chi \in C_{c}^{\infty}(\Omega_2)$ with the property that $\chi \equiv 1$ on some compact neighboorhood $Y' \subseteq Y^{\circ}$ of $Y_0$
and $\rm{supp}(\chi) \subseteq Y$. By Lemma \ref{containlemma}
there exists a neighboorhood $X'$ of $x_0$ and $\varepsilon'\in]0,1]$ such that $f_{\varepsilon}(X') \subseteq Y'$ for $\varepsilon<\varepsilon'$.
Without loss of generality we can assume that $X' \subseteq X$.

In order to show that $\xi_0 \not \in \wf{f^{\ast} u}_{x_0}$ for $x_0 \not\in K_f(u)$ we will show that for a smooth function $\varphi$ with support on a neighboorhood of $x_0$ exists, such that $\mathcal{F}(f^{\ast}u \varphi)$ is rapidly decreasing on $W^c$ (which is a neighboorhood of $\xi_0$).

We choose $\varphi$ to be a smooth function with $\rm{supp}(\varphi) \subseteq X'$. For all $\varepsilon < \varepsilon'$ the
identity $(f_{\varepsilon}^{\ast}) u_{\varepsilon}\ \varphi = f_{\varepsilon}^{\ast}(\chi u_{\varepsilon})\ \varphi$ holds, since the functions $\chi$ and $\varphi$ where choosen
such that $ (\chi \circ f_{\varepsilon})\cdot \varphi \equiv 1$ for $\varepsilon < \varepsilon_1$. 

Set $\widetilde{V}:=\{\eta \in S^{m-1}\mid  \exists \eta_1 \in V, \lambda \in \mathbb{R}_+: \eta=\lambda \cdot \eta_1 \}$ and
$\widetilde{W}:=\{\xi \in S^{m-1}\mid  \exists \xi \in W, \lambda \in \mathbb{R}_+: \xi=\lambda \cdot \xi_1 \}$. 
Obviously $\four{u \chi}$ is rapidly decreasing on $V^c$.
We have that
\begin{eqnarray*}
&&|\four{( f^{\ast} u)_{\varepsilon} \varphi}(\xi)| = |\four{(f^{\ast} \chi u)_{\varepsilon} \varphi } (\xi)| \\
&=&\left| \int_{\mathbb{R}^m} \widehat{\chi u_{\varepsilon}}(\eta) 
\left(\int_{\mathbb{R}^n} \exp{(i \langle f_{\varepsilon} (x), \eta \rangle - i\langle x, \xi \rangle )} \varphi(x) dx \right)d\eta\right|\\
&=&\left|\int_{\mathbb{R}^m} \widehat{\chi u_{\varepsilon}}(\eta) I_{\varepsilon}(\xi, \eta) d\eta \right| 
= \left|\int_{\widetilde{V}} \widehat{\chi u_{\varepsilon}}(\eta) I_{\varepsilon}(\xi, \eta) d\eta \right| + 
\left|\int_{\widetilde{V}^c} \widehat{\chi u_{\varepsilon}}(\eta) I_{\varepsilon}(\xi, \eta) d\eta \right|,
\end{eqnarray*}
where we have set
\begin{eqnarray*}
I_{\varepsilon}(\xi, \eta)  := \int \exp{(i \langle f_{\varepsilon} (x), \eta \rangle - i\langle x, \xi \rangle )} \varphi(x) dx.
\end{eqnarray*}
We intend to apply the stationary phase theorem \ref{statphas} (cf. the Appendix) with
\begin{eqnarray*}
\omega&:=&|\xi|+|\eta| \\
\phi_{\varepsilon} &:=& \langle f_{\varepsilon} (x), \frac{\eta}{|\eta|+|\xi|} \rangle - \langle x, \frac{\xi}{|\eta|+|\xi|} \rangle. 
\end{eqnarray*}
Thus we have to find a bound for the gradient of the phase function $|\phi_{\varepsilon}'(x)| = \left|  {}^T\!df_{\varepsilon}(x) \alpha \eta_1  -  (1-\alpha)\xi_1 \right|$ from below for all $\alpha\in[0,1]$, $\eta_1 = \eta/|\eta|$, $\xi_1=\xi/|\xi|$ and $x \in X$.
By optimization in the parameter $\alpha$ and using the notation of Lemma \ref{notionlemma}  we can bound the expression
\begin{eqnarray*}
\inf_{x \in X,\eta_1 \in V, \xi_1 \in W,\alpha\in[0,1]}  \left| {}^T\!df_{\varepsilon}(x) \alpha \eta_1  -  (1-\alpha) \xi_1 \right|
\end{eqnarray*}
from below by
\begin{eqnarray*}
&& 
\inf_{(x,\eta_1,\xi_1) \in X\times V \times W^c} \frac{| 
{}^{T}\! df_{\varepsilon}(x) \eta_1| }{   |{}^{T}\! df_{\varepsilon}(x) \eta_1 + \xi_1|}
\sqrt{\left|
1- \Big \langle M_{\varepsilon}(x,\eta_1) , \xi_1 
\Big \rangle^2 \right|} \\
&&=
\inf_{(x,\eta_1,\xi_1) \in X\times V \times W^c} \frac{1}{2} \frac{ |  {}^{T}\! df_{\varepsilon}(x) \eta_1| }
{| {}^{T}\! df_{\varepsilon}(x) \eta_1 + \xi_1|}
\left| M_{\varepsilon}(x,\eta_1) -  \xi_1 \right|  \\
&&\ge   \inf_{(x,\eta_1,\xi_1) \in X\times V \times W^c} \frac{1}{4} 
\frac{  |{}^{T}\! df_{\varepsilon}(x) \eta_1| }{
    |{}^T\!df_{\varepsilon}(x) \eta_1| + 1}
\inf_{(x,\eta_1,\xi_1) \in X\times V \times W^c}
  \left|M_{\varepsilon}(x,\eta_1)
 -  \xi_1 \right| \\
 &&\ge   C \sigma_{\varepsilon}^{-1} \inf_{(x,\eta_1,\xi_1) \in X\times V \times W^c}  \left|
M_{\varepsilon}(x,\eta_1)
 -  \xi_1 \right|,  
\end{eqnarray*}
where $\sigma_{\varepsilon}$ is the slow-scaling net from the Definition \ref{unfavdef}. 
According to Lemma \ref{lemmacool} we have that
\begin{equation*} \label{supernet}
d_{\varepsilon} := \inf_{(x,\eta_1,\xi_1) \in X\times V \times W^c}  \left|
M_{\varepsilon}(x,\eta_1)
 -  \xi_1 \right| > d >0 
\end{equation*}
for $\varepsilon< \varepsilon'$, where $d$ is a positive constant.
Thus the gradient of the phase function is uniformly bounded from below by
$|\phi_{\varepsilon}'(x)| \ge  C d \cdot \sigma_{\varepsilon}$ 
for all  $(x,\xi ,\eta) \in X \times \widetilde{V} \times \widetilde{W}^c$. The stationary phase theorem \ref{statphas} yields
\begin{eqnarray}\label{superest1}
\left|I_{\varepsilon}(\xi, \eta) \right| &\le& C_q \varepsilon^{-1} (1+|\xi|+|\eta|)^{-q}
\end{eqnarray}
for all $q \in \mathbb{N}_0$. Note that we use $C, C_p, C_q$ and $C_k,l$ as generic constants.
In the case where $(x, \eta, \xi) \in X \times \widetilde{V}^c \times \widetilde{W}^c$ the stationary phase theorem (now with the phase function $\exp{(-i\langle x, \xi \rangle)}$)
gives
\begin{equation} \begin{split} \label{superexpr}
&|\xi|^{k} |I_{\varepsilon}(\xi, \eta)|=\left|\int  \exp{(- i\langle x, \xi \rangle )}  \exp{(i \langle f_{\varepsilon} (x), \eta \rangle)
 } \varphi(x)  dx \right| \\
&\sum_{|\alpha| \le k} \left|D^{\alpha}_x \left( \exp{(i |\eta| \langle f_{\varepsilon} (x), \eta/|\eta| \rangle )} \varphi \right)\right|  \\ 
\end{split}
\end{equation}
and by repeated use of the chain rule we obtain the estimate
\begin{equation*}
\begin{split}
&\le\sum_{|\alpha| \le k} \sum_{\beta \le \alpha} c_1(\beta) \sup_{x\in X_0} \left|\partial^{\alpha-\beta} \varphi(x)\right| \sum_{l=1}^{|\beta|} |\eta|^l \cdot \\
&\sum_{\gamma_1+...+\gamma_l = \beta} d(\gamma_1, ... ,\gamma_l)  \prod_{1\le i \le l} 
\sup_{x\in X_0} \left|\langle \partial^{\gamma_i} f_{\varepsilon}(x), \eta/|\eta| \rangle\right|
\end{split}
\end{equation*}
where $\gamma_1 + \gamma_2 + \cdots + \gamma_l$ denotes a partition of the multiindex $\beta$ in exactly
$l$ multiindices, that add up componentwise to $\beta$. Using the
notation $|g|_{k} := \max_{|\alpha|=k} |\partial^{\beta} f(x)|$, we can bound the expression by  
\begin{equation*}
\begin{split}
&\le
C_{k,2} \sum_{|\alpha| \le k} \sum_{\beta \le \alpha}  \sum_{l=0}^{|\beta|} |\eta|^l \max_{\gamma_1+...+\gamma_l = \beta}{\left( \prod_{1\le i \le l} 
\sup_{x\in X, \eta \in \widetilde{V}^c} \left|\langle  f_{\varepsilon}(x), \eta/|\eta| \rangle\right|_{|\gamma_i|} \right)} \\
&\le 
C_{k,2} \sum_{|\alpha| \le k} \sum_{\beta \le \alpha}  \sum_{l=0}^{|\beta|} |\eta|^l \max_{0\le j \le |\beta| - l + 1}{    
\left( \sup_{x\in X, \eta \in \widetilde{V}^c} \left|\langle  f_{\varepsilon}(x), \eta/|\eta| \rangle\right|_{j}  \right)^l } \\
&\le
C_{k,3}  (1+|\eta|)^k \max_{l=1}^{k} \max_{0\le j \le k - l + 1}{    
\left( \sup_{x\in X, \eta \in \widetilde{V}^c} \left|\langle  f_{\varepsilon}(x), \eta/|\eta| \rangle\right|_{j}  \right)^l }
.\end{split}
\end{equation*}
Since $f_{\varepsilon}(x)$ is of slow-scale in all derivatives at $x_0$, it holds that for all $j \in \mathbb{N}_0$ there exists 
constants $C_j$ and slow-scaled nets $r_{\varepsilon,j}$ such that
\begin{equation*}
\sup_{(x,\eta_1) \in X\times V^c}{\left|\langle  f_{\varepsilon}(x), \eta_1  \rangle\right|_{j}} \le C_j r_{\varepsilon,j} 
\end{equation*}
holds for small $\varepsilon$.
Summing up the estimates of the form (\ref{superexpr}) give that there exists constants $C_p$ such that
\begin{eqnarray}\label{superest2}
|I_{\varepsilon}(\xi, \eta)| \le C_p (1+|\eta|)^{-p} \varepsilon^{-1} (1+|\xi|)^{p}
\end{eqnarray} 
holds for all $p\in \mathbb{N}_0$ and $(\eta,\xi) \in \widetilde{V}^c \times \widetilde{W}$.  
We observe that
\begin{equation*}\label{last}
\begin{split}
|\four{( f^{\ast} u)_{\varepsilon} \varphi}(\xi)| =& \left|\int_{\mathbb{R}^m} \widehat{\chi u_{\varepsilon}}(\eta) I_{\varepsilon}(\xi, \eta) d\eta \right|\\ 
\le&
C_q  \varepsilon^{-1} \int_{\widetilde{V}} |\widehat{\chi u_{\varepsilon}}(\eta)|  (1+|\xi|+|\eta|)^{-q} d\eta
\\ &+ 
C_p (1+|\xi|)^{-p} \varepsilon^{-1}  \int_{\widetilde{V}^c} |\widehat{\chi u_{\varepsilon}}(\eta)|  (1+|\eta|)^{p} d\eta 
\end{split}
\end{equation*}
holds and using (\ref{superest1}) and (\ref{superest2}) leads to the upper bound
\begin{equation*}
\begin{split}
& C_q \varepsilon^{-1} \sup_{\eta \in \tilde{V} }\left|(1+|\eta|)^{-k} \widehat{\chi u_{\varepsilon}}(\eta) \right|  \int_{\tilde{V}}    (1+|\xi|+|\eta|)^{k-q} d\eta  \\
& + C_p (1+|\xi|)^{-p} \varepsilon^{-1-n} \int_{\widetilde{V}^c} (1+|\eta|)^{p-l}  d\eta. 
\end{split}
\end{equation*}
Finally we set $k:=q-p-m$ and $l:=p+n-1$ and obtain
\begin{equation*}
\begin{split}
& C_q \varepsilon^{-1} \sup_{\eta \in \tilde{V} }\left|(1+|\eta|)^{-k} \widehat{\chi u_{\varepsilon}}(\eta) \right|  \int_{\tilde{V}}    (1+|\xi|+|\eta|)^{k-q} d\eta  \\
& + C_p (1+|\xi|)^{-p} \varepsilon^{-1-n} \int_{\widetilde{V}^c} (1+|\eta|)^{p-l}  d\eta \\
& \le C_{p,q,m,n} \varepsilon^{-1-n} (1+|\xi|)^{-p} 
\end{split}
\end{equation*}
for all $\xi \in \widetilde{W}^c$ and $C_{p,q,m,n}$ some constant depending on $p,q,m$ and $n$.
It follows that $(x_0,\xi_0)  \not \in  \wf{f^{\ast}u}$. 
\end{proof}

\section{Examples}
\begin{example}[Multiplication of Colombeau functions]
This example was presented in \cite[Example 4.2]{HK:01}  in order to show that 
an inclusion relation like in \cite[Theorem 8.2.10]{Hoermander:V1} for the wave front set of a product of distributions,
cannot be extended to Colombeau function with wavefront sets in unfavorable position. 

Consider the Colombeau functions $u$ and $v$ defined by
\begin{eqnarray*}
u_{\varepsilon}&:=& {\varepsilon}^{-1} \rho({\varepsilon}^{-1} (x + \gamma_{\varepsilon} y))\\
v_{\varepsilon}&:=& {\varepsilon}^{-1} \rho({\varepsilon}^{-1} (x - \gamma_{\varepsilon} y)),
\end{eqnarray*}
where $\gamma_{\varepsilon}$ is some net with $\lim_{\varepsilon \rightarrow 0} \gamma_{\varepsilon} = 0$. Note that 
in \cite[Example 4.2]{HK:01} $\gamma_{\varepsilon} := \varepsilon^{1/2}$.
These Colombeau functions are both associated to $\delta(x) \otimes 1(y)$ and the wavefront sets $\wf{u}=\wf{v} = 
\{0\}\times\mathbb{R} \times \{(\pm 1,0)\}$ are in an unfavorable position.
\end{example}
We are going to apply Theorem (\ref{maintheorem}) in order to calculate $WF(u\cdot v)$. First we rewrite $u \cdot v = f^{\ast} \iota(\delta)$,
where $f_{\varepsilon}(x,y)=(x+ \gamma_{\varepsilon} y, x-\gamma_{\varepsilon} y)$
with $\gamma_{\varepsilon}$ some net tending to zero. In Example (\ref{mgexample}) we already showed that $D_f= \mathbb{R}^2 \times S^1 \backslash (\pm 1, \mp1)$.
From \cite[Theorem 15]{HdH:01} it follows that $WF(\iota(\delta)) = \{(0,0)\} \times S^{1}$.
The wavefront unfavorable support of $f$ with respect to $WF(\iota(\delta))$ according to Defintion \ref{unfavdef} is
\begin{eqnarray*}
U_f((0,0)\times S^1) = \{(0,0)\}.
\end{eqnarray*}
Since $f_{\varepsilon}(x,y)$ is of slow-scale in all derivatives at all $x\in\mathbb{R}^2$, it follows that $S_f^c=\emptyset$.
Now Theorem \ref{maintheorem} gives that
\begin{eqnarray*}
\wf{u \cdot v} \subseteq (0,0) \times S^{1}, 
\end{eqnarray*}
which is consistent with the result in \cite[Example 4.2]{HK:01}.

\begin{example}[Hurd-Sattinger]
Let us consider the inital value problem
\begin{equation} \label{ivp} \begin{split}
&\partial_t u + \Theta \partial_x u + \Theta' u  = 0   \\
&u(0,x)  = u_0 \in \col{\mathbb{R}},
\end{split} 
\end{equation}
where $\Theta \in \col{\mathbb{R}^2}$  is defined by $\Theta_{\varepsilon}(x) =  \rho_{\gamma_{\varepsilon}} \ast H(-\cdot)$ with $\rho_{\gamma_{\varepsilon}} = \frac{1}{\gamma_{\varepsilon}}\rho (\cdot /  \gamma_{\varepsilon})$  where $\rho \in \mathcal{S}(\mathbb{R}), \rho\ge 0$ and $\int \rho(x) dx=1$ and $\gamma_{\varepsilon}= \log{(1/\varepsilon)}$ is a net of slow-scale. For the initial value we choose $u_0 := \iota(\delta_{-s_0})$ a delta like singularity at $-s_0$ (for a positive $s_0>0$).
The Hurd-Sattinger example was first given in \cite{HS:68} (it was shown that it is not solvable in $L^1_{\rm{loc}}$, when distributional products are employed). 
It was further investigated in \cite{HdH:01} with methods from Colombeau theory. In \cite{GH:05b} the wavefont set $\rm{WF}_{\gamma}$  (with respect to the slow-scale net $\gamma$) of the Colombeau solution was calculated. For sake of simplicity we do only consider the standard generalized wave-front set $\rm{WF}$, which is smaller since
it neglects the singularities coming from the coefficient $\Theta$.
\end{example}
We have
\begin{eqnarray*}
\Theta_{\varepsilon}(x) = \int_{x/\gamma_{\varepsilon}}^{\infty} \rho(z) dz.
\end{eqnarray*}
Consider the ordinary differential equation for the characteristic curve, which passes through $(t,x)$ at time $t$: 
\begin{eqnarray*}
&&\partial_0 \sigma_{\varepsilon}(s;t,x) = \Theta_{\varepsilon}(\sigma_{\varepsilon}(s;t,x)) \\
&&\sigma(t;t,x)=x.
\end{eqnarray*}
The derivatives of the characteristics fulfill linear initial value problems: If we set $y_{\varepsilon}:=\partial_{t} \sigma_\varepsilon(s;t,x)$
we obtain
\begin{eqnarray*}
\dot{y_{\varepsilon}} = \Theta_{\varepsilon}(\sigma_{\varepsilon}(s;t,x)) y_{\varepsilon},\ y_{\varepsilon}(t;t,x) = -\Theta_{\varepsilon}(x).
\end{eqnarray*}
If we set $z_{\varepsilon}(s;t,x) := \partial_x \sigma_{\varepsilon}(s;t,x)$, we obtain
\begin{eqnarray*}
\dot{z_{\varepsilon}} = -\Theta_{\varepsilon}(\sigma_{\varepsilon}(s;t,x)) z_{\varepsilon},\ z_{\varepsilon}(t;t,x) = 1.
\end{eqnarray*}
This yields
\begin{eqnarray*}
\partial_t \sigma_{\varepsilon}(s;t,x) &=& - \Theta_{\varepsilon}(\sigma_\varepsilon(s;t,x)) \ \text{and} \
\partial_x \sigma_{\varepsilon}(s;t,x) = \frac{\Theta_{\varepsilon}(\sigma_\varepsilon(s;t,x))}{\Theta_\varepsilon(x) }.
\end{eqnarray*}
\newline
If we set $f_{\varepsilon}(t,x):=(t, \sigma_{\varepsilon}(0;t,x))$ then the solution of the initial value problem (\ref{ivp})
reads
\begin{equation*}
u_{\varepsilon} = f_{\varepsilon}^{\ast} (1 \otimes u_{0,\varepsilon})  \cdot \frac{\Theta_\varepsilon(\sigma_\varepsilon(0;t,x))}{\Theta_\varepsilon(x)}.
\end{equation*}
We expect that the wave front set of $u$ is generated from the first factor $f^{\ast}(1 \otimes u_0)$, since
the second factor $\frac{\Theta_\varepsilon(\sigma_\varepsilon(0;t,x))}{\Theta_\varepsilon(x)}$ contains only slow-scale terms.
The wavefront set corresponding to the initial value is $\Gamma:=\rm{WF}(1 \otimes u_0)= \mathbb{R} \times \{-s_0\} \times \{0\} \times \{\pm 1\}$ and
the transposed Jacobian of the 
generalized pullback $f$ is
\begin{eqnarray*}
 df_{\varepsilon}(t,x)^T &&= \left( \begin{array}{cc} 
1 &  - \Theta_{\varepsilon}(\sigma_\varepsilon(0;t,x))\\ 
0 & \frac{\Theta_{\varepsilon} (\sigma_\varepsilon(0;t,x)) } {\Theta_\varepsilon(x) }  
\end{array} \right).
\end{eqnarray*}

If the Jacobian acts on $\zeta := (0, \pm 1) \in S^{1}$ (which are the only irregular directions coming from the wavefront set of $1\otimes u_o$),
 we obtain 
\begin{eqnarray*}
M_{\varepsilon}(t,x,\zeta) &=&\frac{df_{\varepsilon}(t,x)^T \zeta}{|df_{\varepsilon}(t,x)^T \zeta|}\\
&=&\pm \left( \frac{-1}{\sqrt{1+ \Theta_\varepsilon(x)^{-2}}}, \frac{1}{\sqrt{1+\Theta_\varepsilon(x)^2}} \right).
\end{eqnarray*}
Remarkably the result does only depend on the coefficient $\Theta_{\varepsilon}(x)$. Let  $z_{\varepsilon}:=(t_{\varepsilon},x_{\varepsilon}) \rightarrow (t_0,x_0)$ be some convergent net, then we
immediately get that 
\begin{eqnarray*}
M_{\varepsilon}(z_{\varepsilon},\zeta) \rightarrow  \left\{  
\begin{array}{ll} 
(-\zeta_1,\zeta_2) \cdot  &  \text{for}\ x_0<0\\
(0,\zeta_2) \cdot  &  \text{for}\ x_0>0\\ 
\left(-\frac{\alpha}{\sqrt{1+\alpha^2}} \zeta_1,\frac{1}{\sqrt{1+\alpha^2}} \zeta_2 \right)   ,\ \alpha \in [0,1] &\text{for}\ x_0=0.
\end{array} \right.  
\end{eqnarray*}
In order to obtain estimates for the charactistic flow we note that $0 \le \Theta_{\varepsilon} (z) \le 1$,
since the $\rho$ was assumed to be positive. We conclude that
\begin{eqnarray*}
&&\sigma_{\varepsilon}(0;t,x)< \sigma_{\varepsilon}(0;t,y) : x<y \ \text{and} \ \sigma_{\varepsilon}(0;p,x)< \sigma_{\varepsilon}(0;q,x) : q<p, 
\end{eqnarray*}
since $\partial_t \sigma_{\varepsilon}(0;t,x)<0$ and $\partial_x \sigma_{\varepsilon}(0;t,x)>0$.
Thus it is clear that
\begin{eqnarray} \label{charbounds1}
x - t  \le \sigma_{\varepsilon}(0;t,x) \le x
\end{eqnarray}
holds globally. Since $\Theta_{\varepsilon}$ is monotone and positive, we can improve these bounds using the integral representation of $\sigma_{\varepsilon}(0;t,x)$ by
\begin{eqnarray} \label{charbounds}
x- t \Theta_{\varepsilon}(x-t) \le \sigma_{\varepsilon}(0;t,x) \le x - t \Theta_{\varepsilon}(x).
\end{eqnarray}
In order to determine $C_{\Phi_f}$ we have to
calculate all clusterpoints of nets $(\sigma_{\tau(\varepsilon)}(0;t_{\varepsilon},x_{\varepsilon}))_{\varepsilon}$ for arbitrary nets $(t_{\varepsilon},x_{\varepsilon})_{\varepsilon}$ tending to $(t_0,x_0)$ and $\tau:]0,1] \mapsto ]0,1]$ some map with $\tau(\varepsilon)<\varepsilon$.
Assume $t_0>0$, then we have to distinguish three cases:
If $x_0<0$ then we can use the estimate (\ref{charbounds}) to obtain $\lim_{\varepsilon\rightarrow 0} \sigma_{\tau(\varepsilon)}(0;t_{\varepsilon}, x_{\varepsilon})=
x_0-t_0$.
In the case $x_0 > 0$ we first solve $\sigma_{\tau(\varepsilon)}(0;t,x) = x_0 / 2$ 
in the time variable and obtain globally defined (since $\sigma_{\tau(\varepsilon)}(0;t,x)$ is monotone) family of functions $T_{\varepsilon}(x)$ with the property that $\sigma_{\tau(\varepsilon)}(0;T_{\varepsilon}(x), x) = x_0/2$. If $x>x_0/2$ we have that $T_{\varepsilon}(x) >0$.
Without loss of generality we assume that $x_{\varepsilon} > x_0/2$. Using (\ref{charbounds}) we derive the estimate 
\begin{eqnarray*}
T_{\varepsilon}(x) \ge \frac{x - x_{0}/2}{\Theta_{\tau(\varepsilon)}(x_0/2)}.
\end{eqnarray*}
It yields that $\lim_{\varepsilon \rightarrow 0} T_{\varepsilon}(x_{\varepsilon}) = +\infty$, so $t_{\varepsilon} < T_{\varepsilon}(x_{\varepsilon})$
holds for small $\varepsilon$. We obtain
\begin{eqnarray*}
x_{\varepsilon}- t_{\varepsilon} \Theta_{\varepsilon} (x_0/2)\le \sigma_{\varepsilon}(0;t_{\varepsilon},x_{\varepsilon})\le x_{\varepsilon}- t_{\varepsilon} \Theta_{\varepsilon}(x_{\varepsilon})
\end{eqnarray*}
as $\varepsilon$ gets small. It follows that $\lim_{\varepsilon\rightarrow 0} \sigma_{\tau(\varepsilon)}(0;t_{\varepsilon}, x_{\varepsilon}) = x_0$. 

The case $x_0=0$ is a bit more difficult: We proceed with an idea from the proof of \cite[Proposition 3]{O:88}. First off all we can observe that
$-t \le \lim_{\varepsilon \rightarrow 0} \sigma_{\tau(\varepsilon)}(0;t,x) \le 0$.
Let $(a_{\varepsilon})_{\varepsilon}$ be a strictly positive net
with $\lim_{\varepsilon} a_{\varepsilon\rightarrow 0}=0$ and $\lim_{\varepsilon \rightarrow 0} \Theta_{\tau(\varepsilon)} (a_{\varepsilon}) $.
Furthermore we define $t_{\varepsilon}:= t_0 +a_{\varepsilon}$. 
In the next step we construct nets $(x_{\varepsilon,\alpha})_{\varepsilon}$ for each $\alpha\in[-t_0, 0]$ that converge
to zero and have the property that $\lim_{\varepsilon \rightarrow 0} \sigma_{\tau(\varepsilon)}(0;t_\varepsilon, x_{\varepsilon, \alpha})=\alpha$.
We set $x_{\varepsilon, \alpha} := \sigma_{\tau(\varepsilon)}(t_\varepsilon;0,\alpha)$. Let $S_{\varepsilon}(\alpha)$ be the function, that solves
$\sigma_{\tau(\varepsilon)}(S_{\varepsilon}(\alpha);0,\alpha) = a_{\varepsilon}$ globally. It exists because $\partial_{x} \sigma_{\varepsilon}(t;0,x)$ is strictly positive and $\sigma_{\varepsilon}(t;0,x)$ is monotone. Since $a_{\varepsilon} >0$ and
$\alpha \in [-t_0, 0]$ we have that $S_{\varepsilon}(\alpha) \ge 0$. Furthermore we can apply (\ref{charbounds1}) to obtain $S_{\varepsilon}(\alpha) \le -\alpha + a_{\varepsilon}$.
For all $t\ge  S_{\varepsilon}(\alpha)$ it follows that
\begin{equation}
a_{\varepsilon} \le \sigma_{\tau(\varepsilon)}(t;0, \alpha) \le a_{\varepsilon} + (t -S_{\varepsilon}(\alpha))\Theta_{\tau(\varepsilon)}(a_{\varepsilon}) 
\end{equation}
and we can conclude that $\lim_{\varepsilon \rightarrow 0} \sigma_{\tau(\varepsilon)}(t_\varepsilon;0, \alpha) = 0$, since $a_{\varepsilon}$ was choosen such that
$\lim_{\varepsilon \rightarrow0} \Theta_{\tau(\varepsilon)}(a_{\varepsilon}) =0$ and $t_\varepsilon:=t_0 +a_{\varepsilon} \ge  -\alpha + a_{\varepsilon} \ge S_{\varepsilon}(\alpha)$ for $\varepsilon$ small enough.
Finally we have that the generalized graph of the phasefunction $\phi_f$ coming from the characteristic flow $f_{\varepsilon}(t,x)=(t,\sigma_{\varepsilon}(0;t,x))$ is
\begin{eqnarray*}
C_{\phi_f} \cap \tilde{\Gamma} 
&=& \{(t, t-s_0, 0, \pm 1; 0,-s_0, \mp 1, \pm 1)  \mid t<s_0 \} \cup  \\
&&\{ (s_0, 0, 0, \pm 1; 0, -s_0,\mp \frac{\alpha}{\sqrt{1+\alpha^2}}, \pm \frac{1}{\sqrt{1+\alpha^2}} ) \mid \alpha\in[0,1] \} \cup\\
&&\{ (t, 0, 0, \pm 1; t, -s_0, 0, \pm 1 ) \mid t>s_0, \alpha \in [0,1] \} 
\end{eqnarray*}
where $\tilde{\Gamma}:= \{(x,\eta, y, \xi) \in \mathbb{R} \times S^1 \times \mathbb{R} \times S^1 \mid (y,\eta) \in \Gamma \}$.
By Theorem \ref{maintheorem} we obtain
\begin{eqnarray*}
 WF(u) = WF(f^{\ast}(1\otimes u_0))&\subseteq& \pi_2(C_{\phi_f} \cap \tilde{\Gamma}) = \\
&&  \{(t, t-s_0; \mp 1, \pm 1)  \mid t<s_0 \} \cup  \\
&&\{ (s_0, 0; \mp \frac{\alpha}{\sqrt{1+\alpha^2}}, \pm \frac{1}{\sqrt{1+\alpha^2}} ) \mid \alpha\in[0,1] \} \cup\\
&&\{ (t, 0;0, \pm 1 ) \mid t>s_0, \alpha \in [0,1] \}.
 \end{eqnarray*}
\begin{figure} \label{fig}
\caption{Estimated wave front set of the generalized solution of the Hurd-Sattinger partial differential equation.}
\center{\includegraphics[scale=0.30]{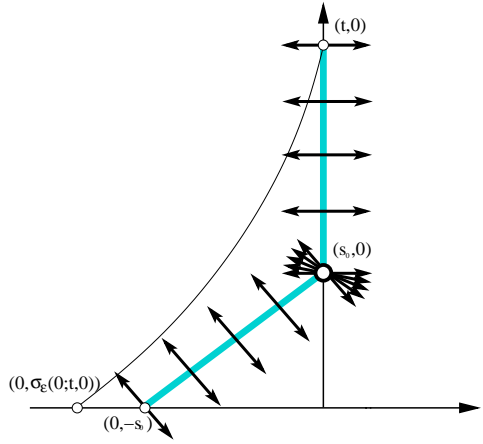}}
\end{figure}

\section{Appendix: Stationary Phase theorem}

\begin{definition} \label{growthnot}
Suppose that $g\in \col{\Omega}$ and let $(g_{\varepsilon})_{\varepsilon}$ be some representative, then we
introduce the notation
\begin{eqnarray*}
|g_{\varepsilon}|_k := \sum_{|\alpha|=k} |\partial^{\alpha} g_{\varepsilon} |
\end{eqnarray*}
and
\begin{eqnarray*}
\mu_{K, k,\varepsilon}(g) &:=& \sup_{K} |g_{\varepsilon}|_k \\
\mu_{K, k,\varepsilon}^{\ast}(g)&:=& \max_{l\le k}{\mu_{K, l, \varepsilon}(g)}.
\end{eqnarray*}
\end{definition}

\begin{lemma} \label{boundderivlemma}
Suppose that $g\in \col{\Omega}$ with an non-negative representative $(g_{\varepsilon})_{\varepsilon}$, such that there exists
$\varepsilon'>0$, with $g_{\varepsilon}(x) \ge 0$ for $\varepsilon <\varepsilon'$ and all $x\in\Omega$. Then it follows
that for any $K \csub \Omega$ some compact set there exist some $C, N$ such that
\begin{eqnarray*}
\delta^2 \sum_{|\alpha|=1} |\partial^{\alpha} g_{\varepsilon}(x)|  
&\le& C \sqrt{g_{\varepsilon}(x)} \sqrt{\mu_{M,2,\varepsilon}^{\ast}(g)}
\end{eqnarray*}
holds for all $x\in K$ and $M$ is some compact set with $K \csub M^{\circ}$.
\end{lemma}
\begin{proof}
Let $K$ be some compact set and $M \csub \Omega$ such that $K \csub M^{\circ}$. We use the notation $B_{\delta}(K):=\{x\in \Omega\mid \exists x_0 \in K: |x-x_0|\le \delta\}$. Then we can find some $\delta >0$ such that $B_{\delta}(K) \subseteq M$ for some $\delta>0$. Then we have for all $x\in K$ that $B_{\delta}(x_0) \subseteq M$.
Choose some $x_0\in K$.
By Taylor's formula we obtain for any $x= x_0+ \alpha v$ with $\alpha \in [0,\delta]$ and $|v|=1$ (note that $x_0+\alpha v \in B_{\delta}(x_0) \subseteq M$), that
\begin{eqnarray*}
0\le g_{\varepsilon}(x) \le g_{\varepsilon}(x_0) + \partial_{v} g_{\varepsilon}(x_0) \alpha + \frac{1}{2} m_\varepsilon \alpha^2,
\end{eqnarray*}
where 
\begin{eqnarray*}
m_{\varepsilon}(v) :=  max( \sup_M |\partial_{v}^2 g_{\varepsilon}|, 2 g_{\varepsilon}(x_0) \delta^{-2} ).
\end{eqnarray*}
So we can obtain the estimate
\begin{eqnarray*}
|\partial_{v} g_{\varepsilon}(x_0)| \le \alpha^{-1} g_{\varepsilon} (x_0) + \frac{\alpha}{2} m_{\varepsilon}(v).
\end{eqnarray*}
Now we distinguish two cases: Let $ \delta^2 m_{\varepsilon} \le 4 g_{\varepsilon}(x_0)$, then we set $\alpha=\delta$ and obtain
\begin{eqnarray*}
\delta|\partial_{v} g_{\varepsilon}(x_0)| &\le& g_{\varepsilon} (x_0) + \frac{\delta^2}{2} m_{\varepsilon}(v) \\
&=&  \sqrt{(g_{\varepsilon} (x_0) + \frac{\delta^2}{2} m_{\varepsilon}(v))^2 } \\
&=&  \sqrt{(g_{\varepsilon}(x_0)^2 + \delta^2 g_{\varepsilon}(x_0) m_{\varepsilon}(v)+ \frac{\delta^4}{4} m_{\varepsilon}(v)^2 } \\
&\le& \sqrt{(g_{\varepsilon}(x_0)^2 +\delta^2 g_{\varepsilon}(x_0) m_{\varepsilon}(v)+ {\delta^2} m_{\varepsilon}(v) g_{\varepsilon}(x_0)} 
\\&=&\sqrt{g_{\varepsilon}(x_0)} \sqrt{ g_{\varepsilon}(x_0) +2 \delta^2 m_{\varepsilon}(v) }.
\end{eqnarray*}
In the case where $ \delta^2 m_{\varepsilon}(v) > 4 g_{\varepsilon}(x_0)$ we set $\alpha:= \sqrt{2 g_{\varepsilon}(x_0)/m_{\varepsilon}(v)}<\delta$
and obtain
\begin{eqnarray*}
|\partial_{v} g_{\varepsilon}(x_0)| &\le& \sqrt{m_{\varepsilon}/2} \sqrt{g_{\varepsilon} (x_0)} + \frac{ \sqrt{2 g_{\varepsilon} (x_0)} m_{\varepsilon}(v) } {2 \sqrt{m_{\varepsilon }(v)} }\\   
&=& \sqrt{2 m_{\varepsilon}(v)}\sqrt{g_{\varepsilon}(x_0)}
\end{eqnarray*}
We can finally conclude that the estimate 
\begin{eqnarray*}
\delta^2 |\partial_{v} g_{\varepsilon}(x_0)|^2
\le g_{\varepsilon}(x_0) ( g_{\varepsilon}(x_0) +2 \delta^2 m_{\varepsilon}(v) ) \le g_{\varepsilon} (x_0) ( \sup_{x\in K} g_{\varepsilon}(x) +2 \delta^2 m_{\varepsilon}(v) )
\end{eqnarray*}
holds for all $x_0\in K$, where $m_{\varepsilon}(v)$ is independant of $x_0$.
\end{proof}

\begin{theorem}[Stationary phase theorem] \label{statphas}
Let $u\in \colcomp{\Omega}$ with support $K \csub \Omega$ and ${\phi}_\varepsilon \in \mathcal{E}(\Omega)$ with the property that
there exists an $\varepsilon_0 >0$ and $m\in \mathbb{N}$ such that
\begin{eqnarray}\label{stphprop}
\inf_{x\in K} (|{\phi}_{\varepsilon}'(x) |) \ge \lambda_{\varepsilon} \ \ 
{\rm for \ all \ }\varepsilon \le \varepsilon_0
\end{eqnarray}
holds, with $\lambda_{\varepsilon}$ some net tending to zero. Then we have that
\begin{eqnarray} \label{statexpr1}
(v_{\varepsilon})_{\varepsilon}:= \left(\int u_{\varepsilon}(x) 
\exp{(i \omega {\phi}_{\varepsilon} (x))} d\!x\right)_{\varepsilon}+\colneg{\Omega}
\end{eqnarray}
is a Colombeau function in the $\omega$ variable and it is bounded by
\begin{eqnarray} \label{statexpr2}
\omega^k |v_{\varepsilon}(\omega)| \le L_{k,\varepsilon} \lambda_{\varepsilon}^{-k} 
\sum_{|\alpha| \le k} \sup_K |D^{\alpha}u_{\varepsilon}| ,
\end{eqnarray}
where  
\begin{eqnarray*}
L_{k,\varepsilon}:=C_k \max{ \{1, \mu^{\ast}_{M,k,\varepsilon}(\phi_{\varepsilon})^{2k^2}   \}}
\end{eqnarray*}
and $M$ is some compact set with $K\subset M^{\circ}$.
\end{theorem}

\begin{proof}
It is obvious that (\ref{statexpr1}) is a well-defined Colombeau function.
The proof follows closely the proof of classical stationary phase theorem in \cite[Theorem 7.7.1]{Hoermander:V1}.
By $N_{\varepsilon}(x):=|\phi_{\varepsilon}'(x)|^2$ we denote the square of the norm of the gradient of the phase function. Let
\begin{eqnarray*}
u_{\nu, \varepsilon} := N^{-1}_{\varepsilon} \frac{\partial \phi_{\varepsilon}}{\partial x_{\nu}} u_{\varepsilon}
\end{eqnarray*}
and since
\begin{eqnarray*}
i \omega \frac{\partial \phi_{\varepsilon}}{\partial x_{\nu}} \exp{(i \omega \phi_{\varepsilon})} = \partial_{\nu} 
\exp{(i \omega \phi_{\varepsilon})}
\end{eqnarray*}
we obtain after an integration by parts 
\begin{eqnarray*}
\int u_{\varepsilon} \exp{(i \omega \phi_{\varepsilon})} dx =
\frac{i}{\omega} \sum_{\nu} \int (\partial_{\nu} u_{\nu, \varepsilon}) \exp{(i \omega \phi_{\varepsilon})} dx
\end{eqnarray*}
(using the notation introduced in Definition \ref{growthnot}).
We prove by induction:
For $k=0$ we have the obvious bound
\begin{eqnarray*}
 |\int u_{\varepsilon} \exp{(i \omega \phi_{\varepsilon})} dx | \le C  \sup_K |u_{\varepsilon}(x)|.
\end{eqnarray*}

Assume that the bound (\ref{statexpr2}) holds for power $k-1$, then we have that
\begin{equation*}\begin{split}
\omega^k | \int u_{\varepsilon} \exp{(i \omega \phi_{\varepsilon})} dx | = &
\omega^{k-1} |\sum_{\nu} \int (\partial_{\nu} u_{\nu, \varepsilon}) \exp{(i \omega \phi_{\varepsilon})} dx|\\
\le& 
L_{k-1,\varepsilon}\sum_{m=0}^{k-1}   
\sup_K  \left(\sum_{\nu=1}^n |u_{\nu, \varepsilon}|_{\mu+1}  N^{m/2-k+1}_{\varepsilon} \right) \label{statexpr8}
\end{split}
\end{equation*}
holds. In the next step we are going to show that 
\begin{eqnarray} \label{statexpr3}
N^{\frac{1}{2}} \sum_{\nu} |u_{\nu, \varepsilon}|_{m}  \le M_{m,\varepsilon} 
\sum_{r=0}^m |u_{\varepsilon}|_{r} N_{\varepsilon}^{\frac{r-m}{2}} 
\end{eqnarray} holds. Again we are prooving by induction: For $m=0$ we have
\begin{eqnarray} \label{statexpr66}
N^{\frac{1}{2}} \sum_{\nu} |u_{\nu, \varepsilon}| = |u_{\varepsilon}| \sum_{\nu}  
N^{-\frac{1}{2}} |\frac{\partial \phi_{\varepsilon} }{\partial x_{\nu} } | \le n |u_{\varepsilon}|
\end{eqnarray}
and let us now assume that (\ref{statexpr3}) holds up to $m-1$. Let $\alpha$ be any multi-index with $|\alpha|=m$.
We apply $\partial^{\alpha}$ on
\begin{eqnarray*}
N_{\varepsilon} u_{\nu, \varepsilon}= u_{\varepsilon} \frac{\partial \phi_{\varepsilon} }{\partial x_{\nu} }
\end{eqnarray*}
and obtain
\begin{eqnarray*}
\partial^{\alpha} (N_{\varepsilon} u_{\nu, \varepsilon})=
\sum_{\beta \le \alpha} {\alpha \choose \beta} (\partial^{\beta} u_{\varepsilon}) 
(\partial^{\alpha-\beta+ e_{\nu}} \phi_{\varepsilon}).
\end{eqnarray*}
It follows that
\begin{eqnarray} \begin{split} \label{statexpr5}
\left|N_{\varepsilon} \partial^{\alpha} u^{\varepsilon}_{\nu}\right| &= 
\left|- \sum_{\beta\ne 0, \beta \le \alpha} {\alpha \choose \beta} \partial^{\beta} N_{\varepsilon}
\partial^{\alpha-\beta} u_{\nu, \varepsilon} + 
\sum_{ \beta \le \alpha} {\alpha \choose \beta} \partial^{\beta} u_{\varepsilon}
\partial^{\alpha-\beta+e_{\nu}} \phi_{\nu}\right|\\
&= C \left( \sum_{l=1}^m |N_{\varepsilon}|_{l} |u_{\nu, \varepsilon}|_{m-l} + 
\sum_{l=0}^{m} |\phi_{\varepsilon}|_{m-l+1} |u_{\varepsilon}|_{l} \right) \\
&= C \left( |N_{\varepsilon}|_{1} |u_{\nu, \varepsilon}|_{m-1} + 
|\phi_{\varepsilon}|_{1} |u_{\varepsilon}|_{m} \right.\\ 
&+ \max_{l=0}^{m-2} {(\sup{|N_{\varepsilon}|_{m-l}}, 
\sup{|\phi_{\varepsilon}|_{m-l+1})}} \sum_{l=0}^{m-2}  (|u_{\nu, \varepsilon}|_{l} 
+  |u_{\varepsilon}|_{l}) \\
&+\left. |u_{\varepsilon}|_{m-1} |\phi_{\varepsilon}|_{2}\right). 
\end{split}
\end{eqnarray} 
Now we can apply Lemma \ref{boundderivlemma} to bound $|N_{\varepsilon}|_1$. This yields 
\begin{eqnarray*} 
|N_{\varepsilon}(x)|_1 &=& \sum_{|\alpha|=1} |\partial^{\alpha} N_{\varepsilon}(x)| 
\le  C_1 \sqrt{N_{\varepsilon}(x)} \sqrt{\mu^{\ast}_{L,2,\varepsilon}(N_{\varepsilon})} \\  
|\phi_{\varepsilon}|_1 &\le& C_2 \sqrt{N_{\varepsilon}(x)}
\end{eqnarray*}
and we can verify that
\begin{eqnarray*}\begin{split}
|\partial^{\alpha} N_{\varepsilon}| &= |\sum_{i=1}^n \partial^{\alpha} (\partial^{e_i} \phi_{\varepsilon}(x))^2|\\
&\le K_1  \sum_{i=1}^n \sum_{\beta \le \alpha} |\partial^{\beta + e_i}  \phi_{\varepsilon}(x)| 
|\partial^{\alpha-\beta + e_i}  \phi_{\varepsilon}(x)| \\
&\le K_2 \sum_{k\le |\alpha|} |\phi_{\varepsilon}(x)|_{|\alpha|-k+1} |\phi_{\varepsilon}(x)|_{k+1}
\end{split}
\end{eqnarray*}
holds. This yields
\begin{equation*}\begin{split}
& \sup_K |N_{\varepsilon}|_l \le K_3 \sum_{k\le l} \sup_K |\phi_{\varepsilon}(x)|_{l-k+1} \sup_K |\phi_{\varepsilon}(x)|_{k+1}\\ 
&\le K_4 \mu^{\ast}_{K,l+1,\varepsilon}(\phi_{\varepsilon})^2.
\end{split}
\end{equation*}
We introduce
\begin{eqnarray*}
\sigma_{\varepsilon,m}&:=&\max_{l=2}^{m} {(\sup{|N_{\varepsilon}|_{l}}, 
\sup{|\phi_{\varepsilon}|_{l+1})}} \\
&\le&\max{(K_4 \mu^{\ast}_{K, m+1,\varepsilon}(\phi_{\varepsilon})^2, \mu^{\ast}_{K, m+1,\varepsilon}(\phi_{\varepsilon}) )}
\end{eqnarray*}
in order to simplify the notation.
Using these bounds we can estimate (\ref{statexpr5}) by
\begin{eqnarray*} \begin{split}
& C \left( |N_{\varepsilon}|^{1/2} C_1 \sqrt{ \mu^{\ast}_{L,2,\varepsilon}(N_{\varepsilon})} |u_{\nu, \varepsilon}|_{m-1} + 
C_2 |N_{\varepsilon}|^{1/2}  
|u_{\varepsilon}|_{m} \right. \\
& \left.+ C_3 \sigma_{\varepsilon}^{(m-2)}  \left( \sum_{l=1}^{m-2}  (|u_{m, \varepsilon}|_{l}  
+  |u_{\varepsilon}|_{l}) \right)+   |\phi_{\varepsilon}|_2 |u_{\varepsilon}|_{m-1} \right) \\
&\le 
C_5 \sigma_{\varepsilon}^{(m-2)}  \left( 
|N_{\varepsilon}|^{1/2} |u_{\nu, \varepsilon}|_{m-1} 
+ |N_{\varepsilon}|^{1/2}  |u_{\varepsilon}|_{m} 
+ \sum_{l=1}^{m-2}  |u_{\nu, \varepsilon}|_{l}     +  \sum_{l=1}^{m-1} |u_{\varepsilon}|_{m-1} \right) 
\end{split}
\end{eqnarray*}
and the induction hypothesis (\ref{statexpr3}) for $m-1$ gives
\begin{eqnarray*}\begin{split} 
&\le 
C_5  \sigma_{\varepsilon}^{(m)}  \cdot \\
&\left( 
  M_{m-1,\varepsilon} \sum_{r=0}^{m-1} |u_{\varepsilon}|_{r} N_{\varepsilon}^{\frac{r-m+1}{2}} 
+ |N_{\varepsilon}|^{1/2}  |u_{\varepsilon}|_{m} \right.
\\& \left.
+ \sum_{l=1}^{m-2}  M_{l,\varepsilon} \sum_{r=0}^l |u_{\varepsilon}|_{r} N_{\varepsilon}^{\frac{r-l-1}{2}} 
  +  \sum_{l=1}^{m-1} |u_{\varepsilon}|_{l} \right) \\
&\le C_5  \sigma_{\varepsilon}^{(m)} \cdot 
\\& \left( 
M_{m-1,\varepsilon} \sum_{r=0}^{m} |u_{\varepsilon}|_{r} N_{\varepsilon}^{\frac{r-m+1}{2}}  
+ \sum_{l=0}^{m-2}  \max{\{1,M_{l,\varepsilon}\}} \sum_{r=0}^{l+1} |u_{\varepsilon}|_{r} N_{\varepsilon}^{\frac{r-l-1}{2}}    
\right) \\ &\le
C_5 \sigma_{\varepsilon}^{(m)}
\left(\sum_{l=0}^{m-1}  \max{\{1,M_{l,\varepsilon}\}} \right) \sum_{r=0}^{l+1} |u_{\varepsilon}|_{r} N_{\varepsilon}^{\frac{r-l-1}{2}} \\
&\le
C_6  \sigma_{\varepsilon}^{(m)} 
\max_{l=0}^{m-1}  \max\{1,M_{l,\varepsilon}\}  \sum_{r=0}^{m} |u_{\varepsilon}|_{r} N_{\varepsilon}^{\frac{r-m-1}{2}}. 
\end{split}\end{eqnarray*}
It follows immediately that 
\begin{eqnarray}
\sum_{\nu} |u_{\nu, \varepsilon}|_{m}  \le M_{m,\varepsilon} \sum_{r=0}^m |u|_{r} N^{\frac{r-m}{2}} 
\end{eqnarray} 
holds, where $M_{m,\varepsilon}$ is defined by
\begin{eqnarray*}
M_{m,\varepsilon} &:=& \sigma_{\varepsilon}^{(m)}  \max_{l=0}^{m-1}  \max{\{1,M_{l,\varepsilon}\}}=C \sigma_{\varepsilon}^{(m)} \Pi_{i=1}^{m-1} \max{\{\sigma_{\varepsilon}^{(i)},1\}}
\end{eqnarray*}
is a generalized number. We are finally able to estimate (\ref{statexpr8}) from above by 
\begin{eqnarray*}
&& L_{k-1,\varepsilon} \sum_{m=0}^{k-1}   
\sup_K  \left(  M_{m+1,\varepsilon} \sum_{r=0}^{m+1} |u_{\varepsilon}|_{r} N_{\varepsilon}^{\frac{r-m-1}{2}}  
N^{\frac{m+1}{2}-k}_{\varepsilon} \right)\\
&\le &
L_{k-1,\varepsilon}M_{k,\varepsilon}  (k-1) \sum_{r=0}^{k}  
\sup_K  \left(   |u_{\varepsilon}|_{r} N_{\varepsilon}^{\frac{r}{2}-k} \right).
\end{eqnarray*}
The generalized constant $L_{k,\varepsilon}$ is recursivly defined by
\begin{eqnarray*}
L_{k,\varepsilon}&:=& L_{k-1,\varepsilon}M_{k,\varepsilon}  (k-1) = (k-1)! \prod_{l=1}^{k} M_{l,\varepsilon}  \\
&=& C (k-1)! \prod_{l=1}^{k}  \sigma_{\varepsilon}^{(l)}  \prod_{i=1}^{l-1} \max{\{\sigma_{\varepsilon}^{(i)},1\}} \\
&\le&
C_k  \left(\prod_{l=1}^{k}  \max{\{1,\mu^{\ast}_{L, l,\varepsilon}(\phi_{\varepsilon})^2\}}    \right) \prod_{l=1}^{k} \prod_{i=1}^{l-1}  \max{\{1, \mu^{\ast}_{L, i,\varepsilon}(\phi_{\varepsilon})^2   \}} \\
&\le& C_k \max{ \{1, \mu^{\ast}_{L, k,\varepsilon}(\phi_{\varepsilon})^{2k^2}   \}} 
\end{eqnarray*}
Note that for $0\le r\le k$ we have that
\begin{eqnarray*}
N_{\varepsilon}^{\frac{r}{2} -k } \le \min{ \{1,\inf_{K} |\phi_{\varepsilon}'(x)|^{-k}\} } \le \lambda_{\varepsilon}^{-k}
\end{eqnarray*}
holds.
\end{proof}

\bibliographystyle{abbrv}
\bibliography{simhal}


\end{document}